\documentclass[aos,seceqn,citesort,MSNbibl,dvips]{arximspdf}

\doi{10.1214/10-AOS815}
\volume{39}
\issue{2}
\pubyear{2011}
\firstpage{931}
\lastpage{955}

\makeatletter

\def\bsuffix #1{#1}

\newtheorem{theorem}{Theorem}
\newtheorem{lemma}{Lemma}[section]
\newproclaim{example}{Example}

\makeatother

\begin{document}
\begin{frontmatter}

\title{A trigonometric approach to quaternary code designs with
application to one-eighth and one-sixteenth fractions}
\runtitle{Quaternary code designs}

\begin{aug}
\author[A]{\fnms{Runchu} \snm{Zhang}\thanksref{t1}\ead[label=e1]{rczhang@nenu.edu.cn}%
\ead[label=e2]{zhrch@nankai.edu.cn}},
\author[B]{\fnms{Frederick K. H.} \snm{Phoa}\thanksref{t2,t3}\ead[label=e3]{fredphoa@stat.sinica.edu.tw}},
\author[C]{\fnms{Rahul} \snm{Mukerjee}\corref{}\thanksref{t1,t4}\ead[label=e4]{rmuk1@hotmail.com}}\\
and
\author[D]{\fnms{Hongquan} \snm{Xu}\thanksref{t2}\ead[label=e5]{hqxu@stat.ucla.edu}}
\runauthor{Zhang, Phoa, Mukerjee and Xu}
\affiliation{Northeast Normal University and Nankai University,
Institute of Statistical Science, Academia Sinica,
Indian Institute of Management Calcutta, and
University of California, Los Angeles}
\address[A]{R. Zhang\\
Klas and School of Mathematics\\
\quad and Statistics\\
Northeast Normal University\\
Changchun 130 024\\
China\\
and\\
LPMC and School\\
\quad of Mathematical Sciences\\
Nankai University\\
Tianjin 300 071\\
China\\
\printead{e1}\\
\phantom{E-mail: }\printead*{e2}} 
\address[B]{F. K. H. Phoa\\
Institute of Statistical Science\hspace*{15.92pt}\\
Academia Sinica\\
Taipei 11529\\
Taiwan\\
\printead{e3}}
\address[C]{R. Mukerjee\\
Indian Institute of Management Calcutta\\
Joka, Diamond Harbour Road\\
Kolkata 700 104\\
India\\
\printead{e4}}
\address[D]{H. Xu\\
Department of Statistics\\
University of California\\
Los Angeles, California 90095-1554\\
USA\\
\printead{e5}}
\end{aug}

\thankstext{t1}{Supported by NNSF China Grants
10871104 and 10571093 and SRFDP China Grant 20050055038.}

\thankstext{t2}{Supported by NSF Grant DMS-08-06137.}

\thankstext{t3}{Supported by National Science Council of Taiwan ROC Grant
98-2118-M-001-028-MY2.}

\thankstext{t4}{Supported by the grant from the Indian
Institute of Management Calcutta.}

\received{\smonth{12} \syear{2009}}

%
\begin{abstract}
The study of good nonregular fractional factorial
designs has received significant attention over the last two decades.
Recent research indicates that designs constructed from quaternary codes
(QC) are very promising in this regard. The present paper shows how a
trigonometric approach can facilitate a systematic understanding of such
QC designs and lead to new theoretical results covering hitherto
unexplored situations. We focus attention on one-eighth and
one-sixteenth fractions of two-level factorials and show that optimal QC
designs often have larger generalized resolution and projectivity than
comparable regular designs. Moreover, some of these designs are found to
have maximum projectivity among all designs.
\end{abstract}

\setattribute{keyword}{AMS}{MSC2010 subject classification.}

\end{frontmatter}

\section{Introduction and preliminaries}\label{sec1}

Fractional factorial designs play a~key role in efficient and economic
experimentation with multiple factors and have gained immense
popularity in various fields of application such as enginee\-ring,
agriculture and medicine. These designs are broadly categorized as
\textit{regu\-lar} and \textit{nonregular} depending on whether or not
they can be generated via defining relations among the factors. In
regular designs, any two factorial effects are either mutually
orthogonal or completely aliased, and the criterion of \textit{maximum
resolution} \cite{1} and its refinement, \textit{minimum aberration}
(MA) \cite{10}, are commonly used in discriminating amongst these
designs. We refer to \cite{14,22} for detailed surveys and extensive
references on regular designs.\looseness=-1

The last two decades, especially the last ten years, have witnessed a
significant spurt in research on nonregular designs. The case of
two-level factors has received particular attention. The notions of
resolution and aberration have been generalized, with statistical
justifications, to these designs; see \cite{7,9,13,19,20,27,28}. More
recent work on nonregular designs and related topics include \cite{4} giving
theoretical results on generalized MA designs, \cite{12} giving a
catalog of
generalized MA designs, \cite{24} on moment aberration projection designs,
\cite{23} on designs obtained from the Nordstrom and Robinson code, and
\cite{17}
on a complete classification of certain two-level orthogonal arrays. It
is well recognized now that although nonregular designs have a~complex
aliasing structure, they can outperform their regular counterparts with
regard to resolution or projectivity, and this is one of the principal
motivating forces for the current surge of interest in these designs.
For more details, see \cite{25} giving a state-of-the-art review of
nonregular designs with a~comprehensive list of references.

A recent major development in nonregular two-level designs has been the
use of quaternary codes (QC) for their efficient construction. The
resulting two-level designs are hereafter called QC designs. While QCs
are known to yield good binary codes in coding theory \cite{11}, QC designs
have been seen to be attractive with regard to resolution, aberration
and projectivity. Moreover, as noted in \cite{25}, these designs are
relatively straightforward to construct and have simple design
representation. Xu and Wong \cite{26} pioneered research on QC designs and
reported theoretical as well as computational results. Phoa and Xu \cite{16}
obtained comprehensive analytical results on quarter fraction QC designs
and showed that they often have larger resolution and projectivity than
regular designs of the same size.

The present paper aims at extending \cite{16} to more highly fractionated
settings. A~serious hurdle in this regard is that the inductive proofs
in \cite{16} become unmanageable when one attempts to go beyond quarter
fractions. A trigonometric representation for QC designs is employed
here in order to overcome this difficulty. This approach is found to be
quite convenient for one-eighth and one-sixteenth fractions which form
our main focus. Earlier, Phoa \cite{15} in his unpublished Ph.D. dissertation
reported partial results on one-sixteenth fraction QC designs under a
certain assumption. An advantage of our approach is that it involves no
such restrictive assumption and enables us to obtain unified and more
comprehensive results on these fractions. As discussed in the concluding
remarks, the trigonometric formulation holds the promise of being
applicable to even more general settings as well.

In Section \ref{sec2}, we introduce the trigonometric approach and present
results pertaining to generalized resolution and wordlength pattern
(WLP) for the case where the number of runs, $N$, is an even power of 2.
The corresponding results when $N$ equals an odd power of 2 appear in
Section \ref{sec3}. Section \ref{sec4} dwells on projectivity and some
directions for
future work are indicated in Section \ref{sec5}. Satisfyingly, with one-eighth
and one-sixteenth fractions, at least over the range covered by our
tables, (a) the same design turns out to be optimal among all QC designs
with respect to all the commonly used criteria like resolution,
aberration and projectivity, and (b) such an optimal design is often
seen to have higher resolution and projectivity than what regular
designs can achieve. The point in (b) reinforces the findings in \cite{16}
for quarter fractions but that in (a) is in contrast to what they
observed in their setup. It is also seen that some of our optimal QC
designs have maximum projectivity among all designs.

Before concluding this section, we reproduce some definitions from \cite{16}
for ease in reference. A two-level design $D$ with $N$ runs and $q$
factors is represented by an $N \times q$ matrix with entries $\pm1$,
where the rows and columns are identified with the runs and factors,
respectively. For any subset $S = \{ c_{1},\ldots,c_{k}\}$ of $k$
columns of~$D$, define
%
%
\begin{equation}\label{equ1.1}
J_{k}(S;D)=\sum_{s = 1}^{N} c_{s1}\cdots c_{sk},
\end{equation}
where $c_{sj}$ is the $s$th entry of column $c_{j}$. The
$J_{k}(S;D)$ values are called the \textit{J-characteristics} of design
$D$; cf. \cite{9,20}. Following \cite{8}, the \textit{aliasing index}
of $S$ is
defined as $\rho_{k}(S;D)=|J_{k}(S;D)|/N$. Clearly, $0 \le\rho_{k}(S;D)
\le1$. If $\rho_{k}(S;D)= 1$, then the columns in $S$ are fully aliased
with one another and form a \textit{complete word} of length $k$ and
aliasing index 1. If $0 <\rho_{k}(S;D)< 1$, then these columns are
partially aliased with one another and form a \textit{partial word} of
length $k$ and aliasing index $\rho_{k}(S;D)$. Finally, if
$\rho_{k}(S;D)= 0$, then these columns do not form a word.

Let $r$ be the smallest integer such that $\max_{\# S =
r}\rho_{r}(S;D)> 0$, where \# denotes cardinality of a set and the
maximum is over all subsets $S$ of $r$ columns of $D$. The
\textit{generalized resolution} \cite{9} of $D$ is defined as
%
%
\begin{equation}\label{equ1.2}
R(D)= r + 1 - \max_{\# S = r}\rho_{r}(S;D).
\end{equation}
For $1 \le k \le q$, let
%
%
\begin{equation}\label{equ1.3}
A_{k}(D)=\sum_{\# S = k} \{ \rho_{k}(S;D)\} ^{2}.
\end{equation}
The vector $(A_{1}(D),\ldots,A_{q}(D))$ is called the \textit{generalized
WLP} of $D$. The \textit{generalized MA} criterion \cite{27}, also
known as
minimum $G_{2}$ aberration \cite{20}, calls for sequential minimization of
$(A_{1}(D), \ldots,A_{q}(D))$. When restricted to regular designs,
generalized resolution, generalized WLP and generalized MA reduce to the
traditional resolution, WLP and MA, respectively. For simplicity, we use
the terminology resolution, WLP and MA for both regular and nonregular
designs.

Following \cite{2}, the design $D$ is said to have \textit
{projectivity} $p$
if every $p$-factor projection contains a complete $2^{p}$ factorial
design, possibly with some points replicated. Evidently, a regular
design of resolution $R$ has projectivity $R - 1$. As shown in \cite{9}, a
design with resolution $R > r$ has projectivity greater than~$r$.

\section{Quaternary code designs in $2^{2n}$ runs}\label{sec2}

\subsection{One-sixteenth fractions}\label{sec21}

In the spirit of \cite{16}, let $C$ be the QC given by the $n \times(n
+ 2)$
generator matrix $[u\enskip v\enskip I_{n}]$, where\vadjust{\goodbreak} $u = (u_{1},\ldots,u_{n})'$ and $v =
(v_{1},\ldots,v_{n})'$ are $n \times1$ vectors over $Z_{4}= \{
0,1,2,3\}$ $(\operatorname{mod} 4)$, $I_{n}$ is the identity matrix of order $n$
over $Z_{4}$, and the primes stand for transpose. The code $C$,
consisting of $4^{n}(\mbox{$=$}2^{2n})$ codewords, each of size $n + 2$, can be
described as
%
%
\begin{equation}\label{equ2.1}
C = \{(a'u,a'v,a_{1},\ldots,a_{n})\dvtx a_{1},\ldots,a_{n}\in Z_{4}\},
\end{equation}
where $a =(a_{1},\ldots,a_{n})'$, and $a'u$ and $a'v$ are reduced mod
4. The
\textit{Gray map}, which replaces each element of $Z_{4}$ with a pair of
two symbols, transforms~$C$ into a binary code $D$, called the binary
image of $C$. For convenience, the two symbols are taken as 1 and $-1$,
instead of the more conventional 0 and~1. Then the Gray map is defined
as
%
%
\begin{equation}\label{equ2.2}\qquad
0\to(1,1),\qquad 1\to(1, -1),\qquad 2\to(-1, -1),\qquad 3\to(-1,1).
\end{equation}
With its codewords as rows, $D$ is a $2^{2n} \times(2n + 4)$ matrix
having entries $\pm1$. Indeed, with columns and rows identified with
factors and runs, respectively, $D$~represents a design involving $2n +
4$ two-level factors and $2^{2n}$ runs. In this sense, $D$ will be
referred to as a $2^{(2n + 4) - 4}$ QC design.

A representation of $D$ using trigonometric functions facilitates the
study of its statistical properties which depend on the choice of $u$
and $v$. Since the pair $(\sqrt{2} \sin(\frac{\pi} {4} + \frac{\pi}
{2}k), \sqrt{2} \cos(\frac{\pi} {4} + \frac{\pi} {2}k))$ equals $(1,1)$,
$(1, -1)$, $(-1, -1)$ and $(-1,1)$ for $k = 0, 1, 2$ and 3 $(\operatorname{mod} 4)$,
respectively, by (\ref{equ2.1}) and (\ref{equ2.2}), the $2^{2n}$ runs in
$D$ can be
expressed as
%
%
\begin{eqnarray}\label{equ2.3}
&&\sqrt{2} \biggl[\sin\biggl(\frac{\pi} {4} + \frac{\pi}
{2}a'u\biggr),\cos\biggl(\frac{\pi} {4} + \frac{\pi} {2}a'u\biggr),\sin\biggl(\frac{\pi} {4}
+ \frac{\pi} {2}a'v\biggr),\nonumber\\
&&\hspace*{21.2pt}\cos\biggl(\frac{\pi} {4} + \frac{\pi} {2}a'v\biggr),
\sin\biggl(\frac{\pi} {4} + \frac{\pi} {2}a_{1}\biggr),\cos\biggl(\frac{\pi} {4} +
\frac{\pi} {2}a_{1}\biggr),\ldots,\nonumber\\[-8pt]\\[-8pt]
&&\hspace*{123.7pt}\sin\biggl(\frac{\pi} {4} + \frac{\pi}
{2}a_{n}\biggr),\cos\biggl(\frac{\pi} {4} + \frac{\pi}
{2}a_{n}\biggr)\biggr],\nonumber\\
&&\eqntext{a_{1},\ldots,a_{n}\in Z_{4}.}
\end{eqnarray}
Denote the $2n + 4$ factors in $D$
by $F_{1},\ldots,F_{4},F_{11},F_{12},\ldots,F_{n1},F_{n2}$, in
conformity with
the ordering in (\ref{equ2.3}), that is, $\sqrt{2} \sin(\frac{\pi} {4} +
\frac{\pi} {2}a'u)$ and $\sqrt{2} \cos(\frac{\pi} {4} + \frac{\pi}
{2}a'u)$ are the levels of $F_{1}$ and $F_{2}$, and so on.

The factors $F_{1},\ldots,F_{4}$, with levels dictated by $u$ or $v$,
require special attention. From this perspective, for any nonempty
collection of factors (or equivalently, columns of $D$), let $x_{k}= 1$
or 0 according as whether $F_{k}$ is included in the collection or not, $1
\le k \le4$. With this notation, the collection is said to be of the
type $x =x_{1}x_{2}x_{3}x_{4}$. Thus, for any binary 4-tuple $x
=x_{1}x_{2}x_{3}x_{4}$, a typical collection of type $x$ consists of
factors $F_{k}$ with $x_{k}= 1$ $(1 \le k \le4)$, and also factors $F_{j1},
F_{j2}$ $(j \in S_{1}), F_{j2}$ $(j \in S_{2})$ and $F_{j1}$ $(j \in S_{3})$,
where\vadjust{\goodbreak} $S_{1},S_{2},S_{3}$ are any disjoint subsets of $\{1,\ldots, n\}$.
The total number of factors in the collection is then $m + X$, where $X
=x_{1} + x_{2} + x_{3} + x_{4}, m = 2n_{1} + n_{2} + n_{3}$ and $n_{j}=
\# S_{j}$ $(1 \le j \le3)$. Here $S_{1},S_{2},S_{3}$ can be empty sets as
well but if $x = 0000$ then at least one of them is nonempty, for
otherwise, the collection contains no factor at all. From (\ref
{equ1.1}), (\ref{equ2.3})
and the definition of aliasing index, it follows that the aliasing index
of a collection of type $x$ $(\mbox{$=$}x_{1}x_{2}x_{3}x_{4})$ as described above
is given by $|V(x)|$, where
%
%
\begin{equation}\label{equ2.4}
V(x)=\sum_{a_{1} = 0}^{3} \cdots\sum_{a_{n} = 0}^{3}
\phi(x;a_{1},\ldots,a_{n})
\end{equation}
with
%
%
\begin{eqnarray}\label{equ2.5}
&&\phi(x;a_{1},\ldots,a_{n})
\nonumber\\
&&\qquad=2^{({1/2})(m + X) - 2n}\sin^{x_{1}}\biggl(\frac{\pi} {4} + \frac{\pi}
{2}a'u\biggr)\cos^{x_{2}}\biggl(\frac{\pi} {4} + \frac{\pi}
{2}a'u\biggr)\nonumber\\[-8pt]\\[-8pt]
&&\qquad\quad{}\times\sin^{x_{3}}\biggl(\frac{\pi} {4} + \frac{\pi}
{2}a'v\biggr)\cos^{x_{4}}\biggl(\frac{\pi} {4} + \frac{\pi} {2}a'v\biggr)
\nonumber\\
&&\qquad\quad{}\times\psi(a_{1},\ldots,a_{n})\nonumber
\end{eqnarray}
and
%
%
\begin{eqnarray}\label{equ2.6}
\psi(a_{1},\ldots,a_{n})
&=& \biggl[\prod_{j \in S_{1}} \biggl\{ \sin\biggl(\frac{\pi} {4} + \frac{\pi}
{2}a_{j}\biggr)\cos\biggl(\frac{\pi} {4} + \frac{\pi} {2}a_{j}\biggr)\biggr\} \biggr]\nonumber\\[-8pt]\\[-8pt]
&&\hspace*{0pt}{}\times\biggl[\prod_{j \in
S_{2}} \cos\biggl(\frac{\pi} {4} + \frac{\pi} {2}a_{j}\biggr)\biggr]\biggl[\prod_{j \in S_{3}}
\sin\biggl(\frac{\pi} {4} + \frac{\pi} {2}a_{j}\biggr)\biggr].\nonumber
\end{eqnarray}
By (\ref{equ2.4})--(\ref{equ2.6}), for any fixed $x$, the value of
$V(x)$ depends on the
sets $S_{1},S_{2},S_{3}$ in addition to $u$ and $v$. Any choice of
$S_{1},S_{2},S_{3}$ that makes $V(x)$ nonzero entails a word of length $m
+ X$ and aliasing index $|V(x)|$. Such a word will be called a word of
type $x$.

We now present Theorem \ref{theo1} below giving an account of words of all
possible types. For $x = 0101$, this result has been proved in the
\hyperref[app]{Appendix}. The proofs for all other $x$ are similar and
occasionally
simpler. In particular, the case $x = 0000$ is evident from the presence
of $I_{n}$ in the generator matrix $[u\enskip v\enskip I_{n}]$ of $C$. Some more
notations will help. With reference to the vectors $u$ and $v$, let
%
%
\begin{eqnarray}\label{equ2.7}\qquad\quad
\Delta_{ks}&=&\{ j\dvtx1 \le j \le n, u_{j} = k, v_{j} = s\},\nonumber\\[-8pt]\\[-8pt]
f_{ks}&=&\#
\Delta_{ks},\qquad 0 \le k,s \le3,
\nonumber\\
\label{equ2.8}
\lambda_{1}&=&f_{10} + f_{30},\qquad \lambda_{2}=f_{01} + f_{03},\qquad
\lambda_{3}=f_{12} + f_{32},\nonumber\\
\lambda_{4}&=&f_{21} + f_{23},\qquad
\lambda_{5}=f_{11} + f_{33},\qquad \lambda_{6}=f_{13} + f_{31},\\
\lambda_{7}&=&f_{02},\qquad
\lambda_{8}=f_{20},\qquad
\lambda_{9}=f_{22},\qquad
\lambda_{10}=f_{00},\nonumber
\\
\label{equ2.9}
l_{1}&=&2(\lambda_{4} + \lambda_{8} + \lambda_{9}) + \lambda_{1} +
\lambda_{3} + \lambda_{5} + \lambda_{6},\nonumber\\
l_{2}&=&2(\lambda_{3} +
\lambda_{7} + \lambda_{9}) + \lambda_{2} + \lambda_{4} + \lambda_{5} +
\lambda_{6},
\nonumber\\
l_{3}&=&2(\lambda_{2} + \lambda_{8} + \lambda_{9}) + \lambda_{1} +
\lambda_{3} + \lambda_{5} + \lambda_{6},\nonumber\\
l_{4}&=&2(\lambda_{1} +
\lambda_{7} + \lambda_{9}) + \lambda_{2} + \lambda_{4} + \lambda_{5} +
\lambda_{6},
\nonumber\\[-8pt]\\[-8pt]
l_{5}&=&2(\lambda_{1} + \lambda_{3} + \lambda_{5} + \lambda_{6}),\qquad
l_{6}=2(\lambda_{2} + \lambda_{4} + \lambda_{5} + \lambda_{6}),\nonumber\\
l_{7}&=&2(\lambda_{1} + \lambda_{2} + \lambda_{3} + \lambda_{4}),
\qquad
l_{8}= 2(\lambda_{7} + \lambda_{8}) + \lambda_{1} + \lambda_{2} +
\lambda_{3} + \lambda_{4},
\nonumber\\
l_{9}&=&2(\lambda_{5} + \lambda_{7} + \lambda_{8}) + \lambda_{1} +
\lambda_{2} + \lambda_{3} + \lambda_{4},\nonumber\\
l_{10}&=&2(\lambda_{6} +
\lambda_{7} + \lambda_{8}) + \lambda_{1} + \lambda_{2} + \lambda_{3} +
\lambda_{4},\nonumber\\
\label{equ2.10}
\rho_{1} &=& 1/2^{\langle({1/2})(\lambda_{1} + \lambda_{3}
+ \lambda_{5} + \lambda_{6})\rangle} ,\qquad \rho_{2} = 1/2^{\langle
({1/2})(\lambda_{2} + \lambda_{4} + \lambda_{5} + \lambda
_{6})\rangle},
\nonumber\\
\xi_{1} &=& 1/2^{\langle({1/2})(\lambda_{1} + \lambda
_{3})\rangle} ,\qquad \xi_{2} = 1/2^{\langle({1/2})(\lambda_{2} +
\lambda_{4})\rangle} ,\\
\xi&=& 1/2^{\langle({1/2})(\lambda_{1} +
\lambda_{2} + \lambda_{3} + \lambda_{4} + 1)\rangle},\nonumber
\end{eqnarray}
where $\langle y\rangle$ is the largest integer not exceeding $y$. By
(\ref{equ2.7}), $f_{ks}$ equals the frequency with which $(ks)$ occurs
as a row of
the $n \times2$ matrix $[u\enskip v]$. Also, each of the quantities introduced
in (\ref{equ2.8})--(\ref{equ2.10}) is uniquely determined by these frequencies.
\begin{theorem}\label{theo1}
With reference to the $2^{(2n + 4) - 4}$ QC
design D, the following hold:

\begin{enumerate}[(a)]
\item[(a)]
For $x = 0000$, there is no word of type $x$.

\item[(b)] For each of $x = 0100$ and 1000, there are
$1/\rho_{1}^{2}$ words of type $x$; every such word has aliasing index
$\rho_{1}$ and length $l_{1} + 1$.

\item[(c)] For each of $x = 0001$ and 0010, there are
$1/\rho_{2}^{2}$ words of type $x$; every such word has aliasing index
$\rho_{2}$ and length $l_{2} + 1$.

\item[(d)] For each of $x = 0111$ and 1011, there are
$1/\rho_{1}^{2}$ words of type $x$; every such word has aliasing index
$\rho_{1}$ and length $l_{3} + 3$.

\item[(e)] For each of $x = 1101$ and 1110, there are
$1/\rho_{2}^{2}$ words of type $x$; every such word has aliasing index
$\rho_{2}$ and length $l_{4} + 3$.

\item[(f)] For $x = 1100$, there is one word of type $x$, with
aliasing index 1 and length $l_{5} + 2$.

\item[(g)] For $x = 0011$, there is one word of type $x$, with
aliasing index 1 and length $l_{6} + 2$.

\item[(h)] For $x = 1111$, there is one word of type $x$, with
aliasing index 1 and length $l_{7} + 4$.

\item[(i)] For each of $x = 0101$ and 1010,
\begin{enumerate}[(i2)]
\item[(i1)] if $\lambda_{5} + \lambda_{6}= 0$, then there are
$1/(\xi_{1}^{2}\xi_{2}^{2})$ words of type $x$; every such word has aliasing
index $\xi_{1}\xi_{2}$ and length $l_{8} + 2$;

\item[(i2)] if $\lambda_{5} + \lambda_{6}> 0$ and $\lambda_{1} +
\lambda_{2} + \lambda_{3} + \lambda_{4}= 0$, then there is one word of
type~$x$, with aliasing index 1 and length $l_{10} + 2$;

\item[(i3)] if $\lambda_{5} + \lambda_{6}> 0$ and $\lambda_{1} +
\lambda_{2} + \lambda_{3} + \lambda_{4}> 0$, then there are
$1/\xi^{2}$ words of type $x$, each with aliasing index $\xi$; half of these
words have length $l_{9} + 2$ and the rest have length $l_{10} +
2$.
\end{enumerate}

\item[(j)] For each of $x = 0110$ and 1001,
\begin{enumerate}[(j2)]
\item[(j1)] if $\lambda_{5} + \lambda_{6}= 0$, then the same conclusion
as in \textup{(i1)} holds;

\item[(j2)] if $\lambda_{5} + \lambda_{6}> 0$ and $\lambda_{1} +
\lambda_{2} + \lambda_{3} + \lambda_{4}= 0$, then there is one word of
type~$x$, with aliasing index 1 and length $l_{9} + 2$;

\item[(j3)] if $\lambda_{5} + \lambda_{6}> 0$ and $\lambda_{1} +
\lambda_{2} + \lambda_{3} + \lambda_{4}> 0$, then the same conclusion as
in \textup{(i3)} holds.
\end{enumerate}
\end{enumerate}
\end{theorem}

Since $\xi=1$ when $\lambda_{1} + \lambda_{2} + \lambda_{3} +
\lambda_{4}= 0$, one can merge \textup{(i2)}, \textup{(i3)}, \textup{(j2)}
and \textup{(j3)} above when
words of all types are considered together. Summarizing Theorem \ref
{theo1}, we
thus get the next result.
\begin{theorem}\label{theo2}
With reference to the $2^{(2n + 4) - 4}$ QC
design D, the following hold:
\begin{enumerate}[(a)]
\item[(a)]
There are $4/\rho_{1}^{2}$ words each with aliasing
index $\rho_{1}$, half of these words have length $l_{1} + 1$ and the
rest have length $l_{3} + 3$.

\item[(b)] There are $4/\rho_{2}^{2}$ words each with aliasing
index $\rho_{2}$; half of these words have length $l_{2} + 1$ and the
rest have length $l_{4} + 3$.

\item[(c)] There are three words each with aliasing index 1;
these have lengths $l_{5} + 2, l_{6} + 2$ and $l_{7} + 4$.

\item[(d)] In addition,
\begin{enumerate}[(d1)]
\item[(d1)] if $\lambda_{5} + \lambda_{6}= 0$, then there are
$4/(\xi_{1}^{2}\xi_{2}^{2})$ words each with aliasing index $\xi_{1}\xi_{2}$
and length $l_{8} + 2$;

\item[(d2)] if $\lambda_{5} + \lambda_{6}> 0$, then there are $4/\xi
^{2}$ words
each with aliasing index $\xi$; half of these words have length $l_{9}
+ 2$
and the rest have length $l_{10} + 2$.
\end{enumerate}
\end{enumerate}
\end{theorem}

Theorem \ref{theo2}, in conjunction with (\ref{equ2.9}) and (\ref
{equ2.10}), shows that the
resolution and WLP of the design $D$ depend on $u$ and $v$ only
through $\lambda_{1},\ldots,\lambda_{10}$. Indeed, for any given $u$ and
$v$, Theorem \ref{theo2} readily yields these features of $D$. This is
illustrated
below.
\begin{example}\label{exam1}
With $n = 3$, let $u = (2,1,1)'$ and $v =(1,1,3)'$. Then
$f_{21} = f_{11} = f_{13} = 1$ and all other $f$'s equal 0, so that by
(\ref{equ2.8}), $\lambda_{4} = \lambda_{5} = \lambda_{6} = 1$ and all
other $\lambda$'s are zeros. Hence by (\ref{equ2.9}) and (\ref{equ2.10}),
\begin{eqnarray*}
l_{1}&=& 4,\qquad l_{2}= 3,\qquad l_{3}= 2,\qquad l_{4}= 3,\qquad l_{5}= 4,\qquad l_{6}= 6,\\
l_{7}&=& 2,\qquad l_{8}= 1,\qquad l_{9}= l_{10}= 3,\qquad \rho_{1}=\rho_{2}=\xi=\tfrac{1}{2}.
\end{eqnarray*}
As a result, parts (a), (b) and (d2) of Theorem \ref{theo2} entail 48
words each
with aliasing index $\frac{1}{2}$; of these, 8 have length four, 32 have
length five and 8 have length 6. Similarly, part (c) of Theorem
\ref{theo2}
entails three words having lengths six, eight, and six, and each with
aliasing index 1. Hence by (\ref{equ1.2}) and (\ref{equ1.3}), in this
case the\vadjust{\goodbreak} QC design~$D$, which is a $2^{10 - 4}$ design, has resolution 4.5 and
WLP $(0,0,0,2,8,4,0,1,0,0)$. As seen later in Table \ref{tab3}, this
design has
maximum resolution and MA among all $2^{10 - 4}$ QC designs. Also, it
has the same WLP but higher resolution than the regular $2^{10 - 4}$ MA
design.
\end{example}

Even though the trigonometric formulation keeps our derivation
tractable, Theorem \ref{theo2} is considerably more involved than its
counterpart,
namely, Theorem 1 of~\cite{16}, for quarter fractions. Consequently,
analytical expressions for the optimal choice of $u$ and $v$, or
equivalently, of $\lambda_{1},\ldots,\lambda_{10}$, maximizing the
resolution or minimizing the aberration of $D$ do not exist in easily
comprehensible forms. On the other hand, as Example \ref{exam1}
demonstrates, for
any given $\lambda_{1},\ldots,\lambda_{10}$, the resolution and WLP of $D$
can be obtained immediately from Theorem \ref{theo2}. Hence, we find
the best
choice of the $\lambda$'s, with regard to resolution and aberration, by
complete enumeration of all possible nonnegative integer-valued
$\lambda_{1},\ldots,\lambda_{10}$ subject to $\lambda_{1} + \cdots+
\lambda_{10} = n$, a condition which is evident from (\ref{equ2.7}) and
(\ref{equ2.8}).
Because of the substantial reduction of the problem as achieved in
Theorem \ref{theo2}, such complete enumeration can be done
instantaneously, for
example by MATLAB, for reasonable values of $n$. The results are
summarized in Table \ref{tab3} and discussed in the next section along
with their
counterparts for QC designs in $2^{2n + 1}$ runs.

\subsection{One-eighth fractions}\label{sec22}

Deletion of any one of the first four columns of the matrix $D$ in
Section \ref{sec21} leads to a QC design involving $2n + 3$ two-level
factors and $2^{2n}$ runs, that is, a $2^{(2n + 3) - 3}$ QC design. The
competing class of QC designs, corresponding to all possible choices of
$u$ and $v$, remains the same up to isomorphism whichever of these four
columns is deleted. This follows by interchanging the roles of $u$ and
$v$ and noting that
\begin{eqnarray*}
\cos\biggl(\frac{\pi} {4} + \frac{\pi} {2}a'u\biggr)&=&\sin\biggl(\frac{\pi} {4} +
\frac{\pi} {2}a'(3u)\biggr),\\
\cos\biggl(\frac{\pi} {4} + \frac{\pi}
{2}a'v\biggr)&=&\sin\biggl(\frac{\pi} {4} + \frac{\pi} {2}a'(3v)\biggr)
\end{eqnarray*}
as $a'u$ and $a'v$ are integers. Therefore, without loss of generality,
we consider the deletion of the first column of $D$. Let $D^{(1)}$ denote
the resulting design. Since the deletion of the first column of $D$
amounts to dropping the factor $F_{1}$, continuing with the notation of
Section \ref{sec21}, only collections of factors of types $x = 0000, 0001,
0010, 0011, 0100, 0101, 0110$ or 0111 can arise now. For all such~$x$,
Theorem \ref{theo1} again describes the numbers of words of type $x$ as
well as
the aliasing indices and lengths of these words. Analogously to
Theorem~\ref{theo2}, this can be summarized as follows.
\begin{theorem}\label{theo3}
With reference to the $2^{(2n + 3) - 3}$ QC
design $D^{(1)}$, the following hold:
\begin{enumerate}[(a)]
\item[(a)]
There are $2/\rho_{1}^{2}$ words each with aliasing
index $\rho_{1}$; half of these words have length $l_{1} + 1$ and the
rest have length $l_{3} + 3$.
\item[(b)] There are $2/\rho_{2}^{2}$ words each with aliasing index
$\rho_{2}$ and length $l_{2} + 1$.\vadjust{\goodbreak}
\item[(c)] There is one word with aliasing index 1 and
length $l_{6} + 2$.
\item[(d)] In addition,
\begin{enumerate}[(d2)]
\item[(d1)] if $\lambda_{5} + \lambda_{6}= 0$, then there are
$2/(\xi_{1}^{2}\xi_{2}^{2})$ words each with aliasing index $\xi_{1}\xi_{2}$
and length $l_{8} + 2$;
\item[(d2)] if $\lambda_{5} + \lambda_{6}> 0$, then there are $2/\xi
^{2}$ words
each with aliasing index $\xi$; half of these words have length $l_{9}
+ 2$
and the rest have length $l_{10} + 2$.
\end{enumerate}
\end{enumerate}
\end{theorem}

Theorem \ref{theo3} shows that the resolution and WLP of $D^{(1)}$
depend on $u$
and $v$ only through $\lambda_{1},\ldots,\lambda_{10}$ and greatly
simplifies the task of finding, by complete enumeration, the
optimal $\lambda$'s maximizing the resolution or minimizing the
aberration of $D^{(1)}$. The results are summarized in Table \ref{tab4} and
discussed in the next section.

\section{Quaternary code designs in 2$^{2n+1}$ runs}\label{sec3}

First, consider one-sixteenth fraction QC designs in $2^{2n + 1}$ runs
as obtained by a branching technique studied in \cite{16} for quarter
fractions. In the present context, this technique can be conveniently
described as follows. Let $\tilde{u}=(u_{0},u_{1},\ldots,u_{n})',
\tilde{v}=(v_{0},v_{1},\ldots,v_{n})'$ be $(n + 1) \times1$ vectors and
$I_{n + 1}$ be the identity matrix of order $n + 1$ over $Z_{4}$.
Consider the QC given by the generator matrix $[\tilde{u}\enskip \tilde{v}\enskip I_{n +
1}]$, and let $\tilde{C}$ be a collection of $2^{2n + 1}$ codewords
thereof, each of size $n + 3$, as given by
%
%
\begin{equation}\label{equ3.1}
\tilde{C}=
\{(\tilde{a}'\tilde{u},\tilde{a}'\tilde{v},a_{0},a_{1},\ldots
,a_{n})\dvtx
a_{0}= 0, 1; a_{1},\ldots,a_{n}\in Z_{4}\},
\end{equation}
where\vspace*{1pt} $\tilde{a} =(a_{0},a_{1},\ldots,a_{n})'$, and
$\tilde{a}'\tilde{u}$ and $\tilde{a}'\tilde{v}$ are reduced mod 4. Apply
the Gray map (\ref{equ2.2}) to $\tilde{C}$ to get a $2^{2n + 1} \times
(2n + 6)$
matrix $\tilde{D}$ having entries $\pm1$. By (\ref{equ3.1}), the
entries in the
fifth and sixth columns of $\tilde{D}$ correspond to the third entry
$a_{0}$ in the codewords of $\tilde{C}$. Since $a_{0}= 0$ or 1, it is
evident from (\ref{equ2.2}) that every entry in the fifth column of
$\tilde{D}$
is 1, while in the sixth column of $\tilde{D}$ half of the entries
equal~1 and the remaining half $-1$. Delete the fifth column of $\tilde{D}$ to
get the final design matrix $D_{0}$, of order $2^{2n + 1} \times(2n +
5)$ and having entries $\pm1$. With its columns and rows identified with
factors and runs, respectively, $D_{0}$ represents a design involving $2n
+ 5$ two-level\vspace*{2pt} factors and $2^{2n + 1}$ runs. In this sense, $D_{0}$ will
be called a $2^{(2n + 5) - 4}$ QC design. Evidently, the role of $u_{0}$
and $v_{0}$ in this construction is different from that of
$u_{1},\ldots,u_{n}$ and $v_{1},\ldots,v_{n}$ and this will be
reflected in
the statistical properties of $D_{0}$. Let $u = (u_{1},\ldots,u_{n})'$
and $v
=(v_{1},\ldots,v_{n})'$.

We consider a trigonometric representation for the runs in $D_{0}$. Since
$D_{0}$ is obtained by deleting the fifth column of $\tilde{D}$, by
(\ref{equ3.1})
and analogously to (\ref{equ2.3}), $D_{0}$ has $2^{2n}$ runs
%
%
\begin{eqnarray}\label{equ3.2}
&&\sqrt{2} \biggl[\sin\biggl(\frac{\pi} {4} + \frac{\pi}
{2}a'u\biggr),\cos\biggl(\frac{\pi} {4} + \frac{\pi} {2}a'u\biggr),\sin\biggl(\frac{\pi} {4}
+ \frac{\pi} {2}a'v\biggr),\nonumber\\
&&\hspace*{21.1pt}\cos\biggl(\frac{\pi} {4} + \frac{\pi}
{2}a'v\biggr),\frac{1}{\sqrt{2}} ,
\sin\biggl(\frac{\pi} {4} + \frac{\pi}
{2}a_{1}\biggr),\nonumber\\[-8pt]\\[-8pt]
&&\hspace*{21.1pt}\cos\biggl(\frac{\pi} {4} +
\frac{\pi} {2}a_{1}\biggr),\ldots,\sin\biggl(\frac{\pi} {4} + \frac{\pi}
{2}a_{n}\biggr),\cos\biggl(\frac{\pi} {4} + \frac{\pi}
{2}a_{n}\biggr)\biggr],\nonumber\\
&&\eqntext{a_{1},\ldots,a_{n}\in Z_{4},}
\end{eqnarray}
which correspond to $a_{0}= 0$, and another $2^{2n}$ runs
%
%
\begin{eqnarray}\label{equ3.3}
&&\sqrt{2} \biggl[\sin\biggl\{ \frac{\pi} {4} + \frac{\pi} {2}(u_{0} + a'u)\biggr\}
,\cos\biggl\{ \frac{\pi} {4} + \frac{\pi} {2}(u_{0} + a'u)\biggr\} ,\nonumber\\
&&\hspace*{20.1pt}\sin\biggl\{
\frac{\pi} {4} + \frac{\pi} {2}(v_{0} + a'v)\biggr\} ,
\cos\biggl\{ \frac{\pi} {4}
+ \frac{\pi} {2}(v_{0} + a'v)\biggr\} ,-
\frac{1}{\sqrt{2}},\nonumber\\[-8pt]\\[-8pt]
&&\hspace*{20.1pt}
\sin\biggl(\frac{\pi} {4} + \frac{\pi} {2}a_{1}\biggr),\cos\biggl(\frac{\pi} {4} +
\frac{\pi} {2}a_{1}\biggr),\ldots,\sin\biggl(\frac{\pi} {4} + \frac{\pi}
{2}a_{n}\biggr),\cos\biggl(\frac{\pi} {4} + \frac{\pi}
{2}a_{n}\biggr)\biggr],\hspace*{-12pt}\nonumber\\
&&\eqntext{a_{1},\ldots,a_{n}\in Z_{4},}
\end{eqnarray}
which correspond to $a_{0}= 1$. Here $a =(a_{1},\ldots,a_{n})'$. Denote
the $2n + 5$ factors in $D_{0}$ by $F_{1},\ldots, F_{4}, F_{5},
F_{11},F_{12},\ldots,F_{n1},F_{n2}$ in conformity with the ordering in
(\ref{equ3.2}) or (\ref{equ3.3}).

In the spirit of Section \ref{sec2}, for any nonempty collection of
factors, let
$x_{k}= 1$ if $F_{k}$ occurs in the collection, and 0 otherwise, $1 \le k
\le5$. The collection is then said to be of the type $(x,x_{5})$, where
$x =x_{1}x_{2}x_{3}x_{4}$. Thus, a typical collection of type
$(x,x_{5})$ consists of factors $F_{k}$ with $x_{k}=1$ $(1 \le k \le5)$,
and also factors $F_{j1},F_{j2}$ $(j \in S_{1}),F_{j2}$ $(j \in S_{2})$ and
$F_{j1}$ $(j \in S_{3})$, where $S_{1},S_{2},S_{3}$ are any disjoint and
possibly empty subsets of $\{1,\ldots, n\}$. Such a collection has $m + X
+ x_{5}$ factors, where $X$ and $m$ are as in Section \ref{sec2}, and
by (\ref{equ1.1}),
(\ref{equ3.2}) and (\ref{equ3.3}), its aliasing index equals $|G(x) + (
- 1)^{x_{5}}H(x)|$,
where
\begin{eqnarray*}
G(x)&=&\frac{1}{2}\sum_{a_{1} = 0}^{3} \cdots\sum_{a_{n} = 0}^{3}
\phi(x;a_{1},\ldots,a_{n}),\\
H(x)&=&\frac{1}{2}\sum_{a_{1} = 0}^{3} \cdots
\sum_{a_{n} = 0}^{3} \phi^{*}(x;a_{1},\ldots,a_{n})
\end{eqnarray*}
with $\phi(x;a_{1},\ldots,a_{n})$ defined by (\ref{equ2.5}) and
$\phi^{*}(x;a_{1},\ldots,a_{n})$ defined similarly replacing $a'u$ and $a'v$
in (\ref{equ2.5}) by $u_{0} + a'u$ and $v_{0} + a'v$. Any choice of
$S_{1},S_{2},S_{3}$ making $G(x) + ( - 1)^{x_{5}}H(x)$ nonzero entails a
word of length $m + X + x_{5}$ and aliasing index $|G(x) + ( -
1)^{x_{5}}H(x)|$. Such a word is called a word of type $(x,x_{5})$.

%
%
\begin{table}
\caption{Values of $N(u_{0}v_{0},\mathit{wl},\mathit{ai})$ for the $2^{(2n + 5) - 4}$ QC
design $D_{0}$}
\label{tab1}
\begin{tabular*}{\tablewidth}{@{\extracolsep{\fill}}lcccccccccc@{}}
\hline
& \multicolumn{10}{c@{}}{$\bolds{u_{0}v_{0}}$}\\[-4pt]
& \multicolumn{10}{c@{}}{\hrulefill}\\
\textbf{\textit{wl},} \textit{\textbf{ai}}
& \textbf{00} & $\bolds{01 / 03}$ & \textbf{02} & $\bolds{10 / 30}$ & $\bolds{11 / 33}$
& $\bolds{12 / 32}$ & $\bolds{13 / 31}$ & \textbf{20} & $\bolds{21
/ 23}$ & \textbf{22}\\
\hline
$l_{1}+1, \theta_{1}$ & 2 & 2 & 2 & 1 & 1 & 1 & 1 & 0 & 0 & 0\\
$l_{1}+2, \theta_{1}$ & 0 & 0 & 0 & 1 & 1 & 1 & 1 & 2 & 2 & 2\\
$l_{2}+1, \theta_{2}$ & 2 & 1 & 0 & 2 & 1 & 0 & 1 & 2 & 1 & 0\\
$l_{2}+2, \theta_{2}$ & 0 & 1 & 2 & 0 & 1 & 2 & 1 & 0 & 1 & 2\\
$l_{3}+3, \theta_{1}$ & 2 & 0 & 2 & 1 & 1 & 1 & 1 & 0 & 2 & 0\\
$l_{3}+4, \theta_{1}$ & 0 & 2 & 0 & 1 & 1 & 1 & 1 & 2 & 0 & 2\\
$l_{4}+3, \theta_{2}$ & 2 & 1 & 0 & 0 & 1 & 2 & 1 & 2 & 1 & 0\\
$l_{4}+4, \theta_{2}$ & 0 & 1 & 2 & 2 & 1 & 0 & 1 & 0 & 1 & 2\\
$l_{5}+2, 1$ & 1 & 1 & 1 & 0 & 0 & 0 & 0 & 1 & 1 & 1\\
$l_{5}+3, 1$ & 0 & 0 & 0 & 1 & 1 & 1 & 1 & 0 & 0 & 0\\
$l_{6}+2, 1$ & 1 & 0 & 1 & 1 & 0 & 1 & 0 & 1 & 0 & 1\\
$l_{6}+3, 1$ & 0 & 1 & 0 & 0 & 1 & 0 & 1 & 0 & 1 & 0\\
$l_{7}+4, 1$ & 1 & 0 & 1 & 0 & 1 & 0 & 1 & 1 & 0 & 1\\
$l_{7}+5, 1$ & 0 & 1 & 0 & 1 & 0 & 1 & 0 & 0 & 1 & 0\\
$l_{8}+2,\omega_{0}$\tabnoteref[*]{t1} & 4 & 2 & 0 & 2 & 0 & 2 & 0 & 0 & 2 & 4\\
$l_{8}+3,\omega_{0}$\tabnoteref[*]{t1} & 0 & 2 & 4 & 2 & 0 & 2 & 0 & 4 & 2 &
0\\[2pt]
$l_{9}+2,\omega$\tabnoteref[\#]{t2} & 2 & 1 & 0 & 1 & 0 & 1 & 2 & 0 & 1 & 2\\
$l_{9}+3,\omega$\tabnoteref[\#]{t2} & 0 & 1 & 2 & 1 & 2 & 1 & 0 & 2 & 1 & 0\\
$l_{10}+2,\omega$\tabnoteref[\#]{t2} & 2 & 1 & 0 & 1 & 2 & 1 & 0 & 0 & 1 & 2\\
$l_{10}+3,\omega$\tabnoteref[\#]{t2} & 0 & 1 & 2 & 1 & 0 & 1 & 2 & 2 & 1 & 0\\
\hline
\end{tabular*}
\tabnotetext[*]{t1}{Entries in these rows, except the ones for $u_{0}v_{0} = 11,
13, 31,
33$, arise if and only if \mbox{$\lambda_{5} + \lambda_{6}= 0$};}
\tabnotetext[\#]{t2}{Entries of these rows, except the ones for $u_{0}v_{0} = 11,
13, 31,
33$, arise if and only if \mbox{$\lambda_{5} + \lambda_{6}> 0$}.}
\end{table}

Steps similar to but more elaborate than those in the \hyperref
[app]{Appendix} may now
be employed to develop an analog of Theorem \ref{theo1} giving an
account of words
of all possible types for the $2^{(2n + 5) - 4}$ QC design $D_{0}$. As
hinted above, this has to be done separately for each possible
pair $u_{0}v_{0}$. One can, thereafter, summarize the findings to get a
counterpart of Theorem \ref{theo2}. However, given the multitude of
possibilities
for the pair $u_{0}v_{0}$, a tabular representation of these summary
results is easier to comprehend than a statement in the form of a
theorem. For any $u_{0}v_{0}$ and any combination of the wordlength
(\textit{wl}) and aliasing index (\textit{ai}), denote the corresponding
number of words in $D_{0}$ by $N(u_{0}v_{0},\mathit{wl},\mathit{ai})/(\mathit{ai})^{2}$.\vadjust{\goodbreak} Table \ref{tab1}
lists all possible $(\mathit{wl},\mathit{ai})$ and, for any such $(\mathit{wl},\mathit{ai})$,
shows $N(u_{0}v_{0},\mathit{wl},\mathit{ai})$ for every $u_{0}v_{0}$. The derivation
underlying this table is omitted to save space. In Table \ref{tab1},
$l_{1},\ldots,l_{10}$ are as in (\ref{equ2.9}), where $\lambda
_{1},\ldots,\lambda_{10}$
continue to be given by (\ref{equ2.8}) with reference to $u =
(u_{1},\ldots,u_{n})'$ and $v =(v_{1},\ldots,v_{n})'$. Also,
%
%
\begin{eqnarray}\label{equ3.4}
\hspace*{32pt}\theta_{1} &=& 1/2^{\langle({1/2})(\lambda_{1} + \lambda
_{3} + \lambda_{5} + \lambda_{6} + \delta_{1})\rangle} ,\qquad \theta_{2} =
1/2^{\langle({1/2})(\lambda_{2} + \lambda_{4} + \lambda_{5} +
\lambda_{6} + \delta_{2})\rangle} ,\nonumber\\
\hspace*{32pt}\omega_{0} &=& \omega_{1}\omega_{2},\qquad
\omega_{1} = 1/2^{\langle({1/2})(\lambda_{1} + \lambda_{3} +
\varepsilon_{1})\rangle} ,\qquad \omega_{2} = 1/2^{\langle
({1/2})(\lambda_{2} + \lambda_{4} + \varepsilon_{2})\rangle} ,\\
\hspace*{32pt}\omega&=& 1/2^{\langle({1/2})(\lambda_{1} + \lambda_{2} + \lambda
_{3} + \lambda_{4} + \varepsilon+ 1)\rangle},\nonumber
\end{eqnarray}
where $\delta_{1} = I$($u_{0}=1$ or 3), $\delta_{2} = I$($v_{0}=1$ or 3),
$\varepsilon_{1} = I$($u_{0}v_{0}=10, 12, 30$ or 32), $\varepsilon_{2} =
I$($u_{0}v_{0}= 01, 03, 21$ or 23), $\varepsilon= \varepsilon_{1} +
\varepsilon_{2}$, and $I(\cdot)$ is the indicator function.

\begin{example}\label{exam2}
With $n = 2$, let $u =(1,2)', v =(2,1)'$ and $u_{0}v_{0}=
11$. Then $f_{12} = f_{21} = 1$ and all other\vadjust{\goodbreak} $f$'s equal 0, so that by
(\ref{equ2.8}), $\lambda_{3} = \lambda_{4} = 1$ and all other $\lambda
$'s are
zeros. Hence, by (\ref{equ2.9}) and (\ref{equ3.4}),
\begin{eqnarray*}
l_{1}&=& l_{2}= 3,\qquad l_{3}= l_{4}= 1,\qquad l_{5}= l_{6}= 2,\qquad l_{7}= 4, \qquad l_{8}=
l_{9}=l_{10}= 2,\\
\theta_{1}&=&\theta_{2}=\tfrac{1}{2},\qquad \omega_{0}=
1,\qquad
\omega=\tfrac{1}{2}.
\end{eqnarray*}
As a result, Table \ref{tab1} shows that there are 48 words each with aliasing
index $\frac{1}{2}$; of these, 24 have length four and 24 have length
five. In addition, there are three words having lengths five, five and
eight, and each with aliasing index 1. Hence by (\ref{equ1.2}) and (\ref
{equ1.3}), in this
case the QC design $D_{0}$, which is a $2^{9 - 4}$ design, has resolution
4.5 and WLP $(0,0,0,6,8,0,0,1,0,0)$. As seen later in Table \ref{tab3}, this
design has maximum resolution and MA among all $2^{9 - 4}$ QC designs.
Also, it has the same WLP but higher resolution than the regular $2^{9 -
4}$ MA design.
\end{example}

We next turn to one-eighth fraction QC designs in $2^{2n + 1}$ runs.
Deletion of any one of the first four columns of $D_{0}$ introduced
earlier leads to a design involving $2n + 4$ two-level factors and
$2^{2n + 1}$ runs, that is, a $2^{(2n + 4) - 3}$ QC design. Using the
same logic as in Section \ref{sec22}, without loss of generality,
%
%
%
\begin{table}[b]
\tabcolsep=0pt
\caption{Values of $N(u_{0}v_{0},\mathit{wl},\mathit{ai})$ for the $2^{(2n + 4) - 3}$
quaternary code design $D_{0}^{(1)}$}
\label{tab2}
\begin{tabular*}{\tablewidth}{@{\extracolsep{\fill}}lcccccccccccccc@{}}
\hline
& \multicolumn{14}{c@{}}{$\bolds{u_{0}v_{0}}$}\\[-4pt]
& \multicolumn{14}{c@{}}{\hrulefill}\\
\textbf{\textit{wl},} \textit{\textbf{ai}}
& \textbf{00} & \multicolumn{1}{c}{$\bolds{01/ 03}$} & \textbf{02} & \textbf{10}
& \textbf{11} & \textbf{12} & \textbf{13} & \textbf{20} &
\multicolumn{1}{c}{$\bolds{21/ 23}$} & \textbf{22} & \textbf{30} & \textbf{31}
& \textbf{32} & \textbf{33}\\
\hline
$l_{1}+1, \theta_{1}$ & 1 & 1 & 1 & $k_{11}$ & $k_{11}$ & $k_{11}$ &
$k_{11}$ & 0 & 0 & 0 & $k_{12}$ & $k_{12}$ & $k_{12}$ & $k_{12}$\\[2pt]
$l_{1}+2, \theta_{1}$ & 0 & 0 & 0 & $k_{12}$ & $k_{12}$ & $k_{12}$ &
$k_{12}$ & 1 & 1 & 1 & $k_{11}$ & $k_{11}$ & $k_{11}$ & $k_{11}$\\[2pt]
$l_{2}+1, \theta_{2}$ & 2 & 1 & 0 & 2 & 1 & 0 & 1 & 2 & 1 & 0 & 2 & 1 &
0 & 1\\[2pt]
$l_{2}+2, \theta_{2}$ & 0 & 1 & 2 & 0 & 1 & 2 & 1 & 0 & 1 & 2 & 0 & 1 &
2 & 1\\[2pt]
$l_{3}+3, \theta_{1}$ & 1 & 0 & 1 & $k_{11}$ & $k_{12}$ & $k_{11}$ &
$k_{12}$ & 0 & 1 & 0 & $k_{12}$ & $k_{11}$ & $k_{12}$ & $k_{11}$\\[2pt]
$l_{3}+4, \theta_{1}$ & 0 & 1 & 0 & $k_{12}$ & $k_{11}$ & $k_{12}$ &
$k_{11}$ & 1 & 0 & 1 & $k_{11}$ & $k_{12}$ & $k_{11}$ & $k_{12}$\\[2pt]
$l_{6}+2, 1$ & 1 & 0 & 1 & 1 & 0 & 1 & 0 & 1 & 0 & 1 & 1 & 0 & 1 & 0\\[2pt]
$l_{6}+3, 1$ & 0 & 1 & 0 & 0 & 1 & 0 & 1 & 0 & 1 & 0 & 0 & 1 & 0 & 1\\[2pt]
$l_{8}+2,\omega_{0}$\tabnoteref[*]{tt1} & 2 & 1 & 0 & $k_{21}$ & 0 & $k_{22}$ & 0 & 0 &
1 & 2 & $k_{22}$ & 0 & $k_{21}$ & 0\\[2pt]
$l_{8}+3,\omega_{0}$\tabnoteref[*]{tt1} & 0 & 1 & 2 & $k_{22}$ & 0 & $k_{21}$ & 0 & 2 &
1 & 0 & $k_{21}$ & 0 & $k_{22}$ & 0\\[2pt]
$l_{9}+2,\omega$\tabnoteref[\#]{tt2} & 1 & $\frac{1}{2}$ & 0 & $\frac{1}{2}$ & 0 &
$\frac{1}{2}$ & 1 & 0 & $\frac{1}{2}$ & 1 & $\frac{1}{2}$ & 1 &
$\frac{1}{2}$ & 0\\[2pt]
$l_{9}+3,\omega$\tabnoteref[\#]{tt2} & 0 & $\frac{1}{2}$ & 1 & $\frac{1}{2}$ & 1 &
$\frac{1}{2}$ & 0 & 1 & $\frac{1}{2}$ & 0 & $\frac{1}{2}$ & 0 &
$\frac{1}{2}$ & 1\\[2pt]
$l_{10}+2,\omega$\tabnoteref[\#]{tt2} & 1 & $\frac{1}{2}$ & 0 & $\frac{1}{2}$ & 1 &
$\frac{1}{2}$ & 0 & 0 & $\frac{1}{2}$ & 1 & $\frac{1}{2}$ & 0 &
$\frac{1}{2}$ & 1\\[2pt]
$l_{10}+3,\omega$\tabnoteref[\#]{tt2} & 0 & $\frac{1}{2}$ & 1 & $\frac{1}{2}$ & 0 &
$\frac{1}{2}$ & 1 & 1 & $\frac{1}{2}$ & 0 & $\frac{1}{2}$ & 1 &
$\frac{1}{2}$ & 0\\[2pt]
\hline
\end{tabular*}
\tabnotetext[*]{tt1}{Entries in these rows, except the ones for $u_{0}v_{0} = 11, 13, 31, 33$,
arise if and only if \mbox{$\lambda_{5} + \lambda_{6}= 0$};}
\tabnotetext[\#]{tt2}{Entries of these rows, except the ones for $u_{0}v_{0} = 11, 13, 31, 33$,
arise if and only if \mbox{$\lambda_{5} + \lambda_{6}> 0$}.}
\end{table}
suppose the
first column of $D_{0}$ is deleted. Let $D_{0}^{(1)}$ denote the
resulting design. Table \ref{tab2} lists all possible combinations
$(\mathit{wl},\mathit{ai})$ of
the wordlength and aliasing index in $D_{0}^{(1)}$ and, for any
such $(\mathit{wl},\mathit{ai})$, shows $N(u_{0}v_{0},\mathit{wl},\mathit{ai})$ for every $u_{0}v_{0}$, where
$N(u_{0}v_{0},\mathit{wl},\mathit{ai})$ is defined as above but now refers to
$D_{0}^{(1)}$. In Table \ref{tab2}, $l_{1},l_{2}$ etc. are as in (\ref{equ2.9}),
$\theta_{1},\theta_{2},\omega_{0}$ and $\omega$ are as in (\ref
{equ3.4}), and
\begin{eqnarray*}
k_{11}&=&\tfrac{1}{2}I(\lambda_{1} + \lambda_{3} > 0),\qquad k_{12} = 1 -
k_{11},\\
k_{21} &=& I(\lambda_{1} + \lambda_{3} + \lambda_{5} + \lambda_{6}
> 0),\qquad k_{22} = 2 - k_{21}.
\end{eqnarray*}

As illustrated in Example \ref{exam2}, Tables \ref{tab1} and \ref{tab2}
readily yield, in their
respective contexts, the resolution and WLP of a QC design for any given
$\lambda_{1},\ldots,\lambda_{10}$ and $u_{0}v_{0}$. Hence, we find the best
choice of the $\lambda$'s and $u_{0}v_{0}$, with regard to resolution and
aberration, by complete enumeration of all possibilities. Again, because
of the significant reduction achieved in Tables \ref{tab1} and \ref
{tab2} by theoretical
means, such complete enumeration can be done instantaneously, for
example by MATLAB, for reasonable values of $n$. The results are
summarized in Tables \ref{tab3} and \ref{tab4} for one-sixteenth and
one-eighth fractions,
respectively.

%
%
\begin{table}
\tabcolsep=0pt
\caption{One-sixteenth fraction QC designs with maximum resolution and
MA}
\label{tab3}
\begin{tabular*}{\tablewidth}{@{\extracolsep{\fill}}lll@{}}
\hline
& & \multicolumn{1}{c@{}}{\textbf{Regular}}\\
\textbf{Design} & \multicolumn{1}{c}{\textbf{QC design with maximum resolution and MA}} & \multicolumn{1}{c@{}}{\textbf{MA design}} \\
\hline
$2^{8-4}$ & $\lambda= 0011000000$, $R = 4$, $A = (14, 0, 0, 0, 1)$ & $R =
4$, $A$ same\\
$2^{9 - 4}$ & $\lambda= 0011000000$, $u_{0}v_{0}= 11$, $R = 4.5$, $A = (6,
8, 0, 0, 1, 0)$ & $R = 4$, $A$ same\\
$2^{10 - 4}$ & $\lambda= 0001110000$, $R= 4.5$, $A = (2, 8, 4, 0, 1, 0,
0)$ & $R = 4$, $A$ same\\
$2^{11 - 4}$ & $\lambda= 0001110000$, $u_{0}v_{0}= 12$, $R = 5.5$, $A = (0,
6, 6, 2, 1, 0, 0, 0)$ & $R = 5$, $A$ same\\
$2^{12 - 4}$ & $\lambda= 0011110000$, $R = 6.5$, $A = (0, 0, 12, 0, 3, 0,
0, 0, 0)$ & $R = 6$, $A$ same\\
$2^{13 - 4}$ & $\lambda= 0011110000$, $u_{0}v_{0}= 22$, $R = 6.5$, $A = (0,
0, 4, 8, 3, 0, 0, 0, 0, 0)$ & $R = 6$, $A$ same\\
$2^{14 - 4}$ & $\lambda= 1011110000$, $R= 6.5$, $A = (0, 0, 2, 8, 3, 0, 2,
0, 0, 0, 0)$ & $R = 7$, $A$ better\\
\hline
\end{tabular*}
\end{table}

%
%
\begin{table}[b]
\tabcolsep=0pt
\caption{One-eighth fraction QC designs with maximum resolution and MA}
\label{tab4}
\begin{tabular*}{\tablewidth}{@{\extracolsep{\fill}}lll@{}}
\hline
& & \multicolumn{1}{c@{}}{\textbf{Regular}}\\
\textbf{Design} & \multicolumn{1}{c}{\textbf{QC design with maximum resolution and MA}}
& \multicolumn{1}{c@{}}{\textbf{MA design}}\\
\hline
$2^{7 - 3}$ & $\lambda= 0011000000$, $R = 4$, $A = (7, 0, 0, 0 )$ & $R =
4$, $A$ same\\
$2^{8 - 3}$ & $\lambda= 0011000000$, $u_{0}v_{0}= 11$, $R = 4.5$, $A = (3,
4, 0, 0, 0)$ & $R = 4$, $A$ same\\
$2^{9 - 3}$ & $\lambda= 0010110000$, $R= 4.5$, $A = ( 1, 4, 2, 0, 0, 0)$ &
$R = 4$, $A$ same\\
$2^{10 - 3}$ & $\lambda= 0010110000$, $u_{0}v_{0}= 21$, $R= 5.5$, $A = (0,
3, 3, 1, 0, 0, 0)$ & $R = 5$, $A$ same\\
$2^{11 - 3}$ & $\lambda= 0011110000$, $R= 6.5$, $A = (0, 0, 6, 0, 1, 0, 0,
0)$ & $R = 6$, $A$ same\\
$2^{12 - 3}$ & $\lambda= 0011110000$, $u_{0}v_{0}= 12$, $R= 6.75$, $A = (0,
0, 2, 4, 1, 0, 0, 0, 0)$ & $R = 6$, $A$ same\\
$2^{13 - 3}$ & $\lambda= 0021110000$, $R= 7.75$, $A = ( 0, 0, 0, 4, 3, 0,
0, 0, 0, 0)$ & $R = 7$, $A$ same\\
\hline
\end{tabular*}
\end{table}

A brief discussion on Tables \ref{tab3} and \ref{tab4}, showing
one-sixteenth and
one-eighth fraction QC designs with maximum resolution and MA is in
order. The 14 QC designs shown in these tables are optimal, among all
comparable QC designs, under both criteria. All these 14 designs have
resolution four or higher, and for each, the resolution $R$ and $A =
(A_{4},A_{5},\ldots)$ are shown. For ease in comparison, we also show $R$
and comment on $A$ for the corresponding regular MA designs, as obtained
by Chen and Wu \cite{6}. Out of the 14 optimal QC designs in Tables
\ref{tab3}
and \ref{tab4},
there are two, that is, the first design in either table, which have the
same $R$ and $A$ as the corresponding regular MA designs. It can be seen
that these two designs involve only full words and hence are themselves
regular. A comparison with Table 6 of Sun, Li and Ye \cite{18} shows that
these two designs have maximum resolution and MA in their sense among
all designs of the same size. Eleven of the remaining twelve optimal QC
designs in our Tables \ref{tab3} and \ref{tab4} have the same WLP but
higher resolution
than the corresponding regular MA designs. Only the $2^{14 - 4}$ optimal
QC design turns out to be worse than the regular MA design. In both
Tables \ref{tab3} and \ref{tab4}, $\lambda$ stands for the
10-tuple $\lambda_{1}\lambda_{2}\cdots\lambda_{10}$.

Indeed, the theoretical results reported in Theorems \ref{theo2}, \ref
{theo3} and Tables~\ref{tab1},~\ref{tab2} readily allow extension of Tables \ref{tab3} and \ref{tab4}
beyond the ranges considered
here, if the situation so demands. For instance, from Table \ref{tab2},
one can
check that the $2^{16 - 3}$ QC design with maximum resolution and MA is
given by $\lambda= 0020220000, u_{0}v_{0}= 20$. This design has $R =
8.875$ and $A = (0, 0, 0, 0, 1, 4, 2, 0,\break 0, 0, 0, 0, 0)$, while the
corresponding regular MA design has the same $A$ but $R = 8$.

\section{Results on projectivity}\label{sec4}

The following results give upper
bounds on the projectivity of the one-sixteenth fraction QC designs $D$
and $D_{0}$ introduced in Sections \ref{sec2} and \ref{sec3}, respectively.
\begin{theorem}\label{theo4}
The projectivity $p$ of the $2^{(2n + 4) -
4}$ QC design D satisfies \textup{(i)} $p \le\frac{4}{3}n + 1$, if $n
= 0
\operatorname{mod} 3$, \textup{(ii)} $p \le\frac{4}{3}(n - j) + 3$,
if $n = j \operatorname{mod} 3$, with $j = 1$ or
$2$.
\end{theorem}
\begin{pf}
We prove only (ii). The proof of (i) is similar. Let $n =
3t + j$, where $t$ is an integer and $j = 1$ or 2. We need to show
that $p \le4t + 3$. If $p \ge4t + 4$, then all full words in $D$ have
length at least $4t + 5$, so that by Theorem \ref{theo2}(c) and (\ref{equ2.9}),
\begin{eqnarray*}
l_{5}+2 &=& 2(\lambda_{1} + \lambda_{3} + \lambda_{5} + \lambda_{6}) +
2\ge4t + 5,\\
l_{6}+2 &=& 2(\lambda_{2} + \lambda_{4} + \lambda_{5} +
\lambda_{6}) + 2\ge 4t+5,
\\
l_{7}+4 &=& 2(\lambda_{1} + \lambda_{2} + \lambda_{3} + \lambda_{4}) +
4\ge 4t+5,
\end{eqnarray*}
that is, invoking the integrality of $\lambda_{1},\ldots,\lambda_{6}$,
\begin{eqnarray*}
\lambda_{1} + \lambda_{3} + \lambda_{5} + \lambda_{6} &\ge&2t +
2,\qquad
\lambda_{2} + \lambda_{4} + \lambda_{5} + \lambda_{6} \ge 2t + 2,\\
\lambda_{1} + \lambda_{2} + \lambda_{3} + \lambda_{4} &\ge&2t + 1.
\end{eqnarray*}
Adding the last three inequalities, $2(\lambda_{1} + \cdots+ \lambda_{6})
\ge6t + 5$, that is, $n \ge\lambda_{1} + \cdots+ \lambda_{6} \ge3t + 3$,
again using the integrality of $\lambda_{1},\ldots,\lambda_{6}$, and we
reach a contradiction.
\end{pf}

From Table \ref{tab1}, now observe that the $2^{(2n + 5) - 4}$ QC
design $D_{0}$
has at least three full words. These have lengths (a) $l_{5}+2$, $l_{6}+2,
l_{7}+4$ if $u_{0}v_{0}= 00, 02, 20, 22$, (b) $l_{5}+2$, $l_{6}+3$,
$l_{7}+5$ if $u_{0}v_{0}= 01, 03, 21, 23$, (c) $l_{5}+3$, $l_{6}+2$,
$l_{7}+5$ if $u_{0}v_{0} = 10, 12, 30, 32$ and (d) $l_{5}+3$, $l_{6}+3$,
$l_{7}+4$ if $u_{0}v_{0}= 11, 13, 31, 33$. Hence, arguments similar to
but more elaborate than those in Theorem \ref{theo4} lead to the
following result.
\begin{theorem}\label{theo5}
The projectivity p of the $2^{(2n + 5) - 4}$ QC
design $D_{0}$ satisfies \textup{(i)} $p \le\frac{4}{3}n + 2$ if $n =
0 \operatorname{mod} 3$, \textup{(ii)} $p \le\frac{4}{3}(n - 1) + 3$,
if $n = 1 \operatorname{mod} 3$,
\textup{(iii)} $p \le\frac{4}{3}(n - 2) + 4$, if $n = 2 \operatorname
{mod} 3$.
\end{theorem}

Table \ref{tab5} shows the projectivities of the one-sixteenth
fraction QC
designs reported in Table \ref{tab3}. It is easily seen that these
designs attain
the upper bounds on projectivity as shown in Theorems \ref{theo4} or
\ref{theo5}. Hence, in
addition to having maximum resolution and MA, they have maximum
projectivity among all comparable QC designs. Indeed, for $8 \le q \le
12$, the $2^{q - 4}$ QC designs in Table \ref{tab3} have projectivity
$q -
5$ which is the highest among all designs of the same size. This holds
because, otherwise, one would get a $2^{q - 4}$ design with
projectivity $q - 4$, that is, an orthogonal array OA $(2^{q - 4}, q, 2,q
- 4)$ of index unity, which is nonexistent; see~\cite{3}. Table~\ref
{tab5} also shows
that the use of QC designs leads to gain in projectivity over regular MA
designs for $9 \le q \le14$.

%
%
\begin{table}[b]
\caption{Projectivities of the one-sixteenth fraction QC designs in Table
\protect\ref{tab3}}
\label{tab5}
\begin{tabular*}{\tablewidth}{@{\extracolsep{\fill}}lccccccc@{}}
\hline
\textbf{Design} & \multicolumn{1}{c}{$\bolds{2^{8 - 4}}$} & $\bolds{2^{9 - 4}}$
& $\bolds{2^{10 - 4}}$ & $\bolds{2^{11 - 4}}$ &
$\bolds{2^{12 - 4}}$ & $\bolds{2^{13 - 4}}$ & $\bolds{2^{14 - 4}}$\\
\hline
Projectivity of the & 3 & 4 & 5 & 6 & 7 & 7 & 7 \\
\quad QC design in Table \ref{tab3} \\
Projectivity of the & 3 & 3 & 3 & 4 & 5 & 5 & 6 \\
\quad regular MA design \\
\hline
\end{tabular*}
\end{table}

We next consider one-eighth fractions and show in Table \ref{tab6} the
projectivities of the QC designs reported in Table \ref{tab4}. For $7
\le q \le
11$, the $2^{q - 3}$ QC designs in Table~\ref{tab4} are seen to have
projectivity
$q - 4$ which is the highest among all designs of the same size. This
follows as in the last paragraph using a nonexistence result in \cite
{3} on
orthogonal arrays of index unity. Also, for $q = 12$ and 13, the QC
designs in Table~\ref{tab4} were computationally verified to have maximum
projectivity at least among all comparable QC designs. Incidentally, for
one-eighth fraction QC designs, it is hard to develop analogs of
Theorems \ref{theo4} and \ref{theo5} as there is only one guaranteed
full word (cf. Theorem
\ref{theo3} and Table \ref{tab2}) but the computational study of
projectivity remains
manageable with a moderate number of factors. Table \ref{tab6} also
shows the
projectivity of regular MA designs and the gains via the use of QC
designs, for $8 \le q \le13$, are evident.

%
%
\begin{table}
\caption{Projectivities of the one-eighth fraction QC designs in Table
\protect\ref{tab4}}
\label{tab6}
\begin{tabular*}{\tablewidth}{@{\extracolsep{\fill}}lccccccc@{}}
\hline
\textbf{Design} & $\bolds{2^{7 - 3}}$ & $\bolds{2^{8 - 3}}$ & $\bolds{2^{9 - 3}}$ & $\bolds{2^{10 - 3}}$ &
$\bolds{2^{11 - 3}}$ & $\bolds{2^{12 - 3}}$ & $\bolds{2^{13 - 3}}$\\
\hline
Projectivity of the & 3 & 4 & 5 & 6 & 7 & 7 & 7\\
\quad QC design in Table \ref{tab4} \\
Projectivity of the & 3 & 3 & 3 & 4 & 5 & 5 & 6\\
\quad regular MA design \\
\hline
\end{tabular*}
\end{table}

The foregoing discussion reveals that, unlike what often happens with
quarter fraction QC designs \cite{16}, maximum projectivity is not in
conflict with maximum resolution or MA in our setup at least over the
range covered by Tables \ref{tab3}--\ref{tab6}, where the same QC
design turns out to be
optimal, among all comparable QC designs, with regard to all the three
criteria.

\section{Summary and future work}\label{sec5}

In the present paper, a trigonometric approach was developed to obtain
theoretical results on QC designs with focus on one-eighth and
one-sixteenth fractions of two-level factorials. It was seen that
optimal QC designs often have larger resolution and projectivity than
comparable regular designs. In addition, some of these designs were
found to have maximum projectivity among all designs.

Before concluding, we indicate a few open issues. It should be possible
to use the trigonometric approach to obtain further theoretical results
on the projectivity of QC designs. This calls for examining the
existence of solutions to certain trigonometric equations. For instance,
by (\ref{equ2.3}), the first four columns of the $2^{(2n + 4) - 4}$ QC
design $D$
in Section \ref{sec2} contain a full 2$^{4}$ factorial if and only if for
every $y_{1},\ldots,y_{4}$ in $\{ - 1,1\}$, the equations
\begin{eqnarray*}
\sqrt{2} \sin\biggl(\frac{\pi} {4} + \frac{\pi} {2}a'u\biggr) &=& y_{1},\qquad \sqrt{2}
\cos\biggl(\frac{\pi} {4} + \frac{\pi} {2}a'u\biggr) = y_{2},
\\
\sqrt{2} \sin\biggl(\frac{\pi} {4} + \frac{\pi} {2}a'v\biggr) &=& y_{3},\qquad \sqrt{2}
\cos\biggl(\frac{\pi} {4} + \frac{\pi} {2}a'v\biggr) = y_{4},
\end{eqnarray*}
admit a solution for $a = (a_{1},\ldots,a_{n})'$ in $Z_{4}$. A study of
equations of this kind, however, branches out into too many cases,
depending on $u$ and $v$, compared to the derivation of results on
wordlength and aliasing index as done here.

It would also be of interest to investigate how the trigonometric
approach can be implemented for even more highly fractionated QC
designs. The trigonometric formulation as well as the mathematical tools
are expected to be essentially same as the ones here and the main
difficulty will lie in handling the multitude of cases that such an
effort will involve. A related issue concerns the development of a
complementary design theory for QC designs in the spirit of \cite{5,21} and
with respect to an appropriately defined reference set.

Some kind of symbolic computation may help in addressing the open
problems mentioned above. We hope that the present endeavor will
generate further interest in these and related issues.

\begin{appendix}\label{app}

\section*{Appendix: Proof of Theorem \lowercase{\protect\ref{theo1}} for $x = 0101$}

The proof will be worked out through a sequence of lemmas. We
concentrate on $V(x)$ and begin by giving an expression for
$\psi(a_{1},\ldots,a_{n})$ in (\ref{equ2.6}). Recall that in (\ref{equ2.6}),
$S_{1},S_{2},S_{3}$ are disjoint subsets of $\{1,\ldots, n\}$ and that $m
=2n_{1} + n_{2} + n_{3}$ with $n_{k}=\# S_{k}$. For $1 \le k \le3$, let
$\Sigma_{k}$ denote the sum over $2^{n_{k}}$ terms corresponding to the
$2^{n_{k}}$ subsets $W_{k}$ of $S_{k}$, and for any such subset $W_{k}$,
write $\bar{W}_{k}=S_{k}\setminus W_{k},w_{k} = \# W_{k} ,\bar{w}_{k} =
\# \bar{W}_{k}$. Thus if $S_{1}= \{2, 3\}$, then $\Sigma_{1}$ denotes
the sum over 2$^{2}$ terms corresponding to $W_{1}= \mbox{empty set}$,
$\{2\}$,
$\{3\}$ and $\{2, 3\}$. For any given subsets $W_{1},W_{2},W_{3}$
of $S_{1},S_{2},S_{3}$, we also write $\Sigma^{(1)}, \bar{\Sigma}
^{(1)},\Sigma^{(23)}$ and $\bar{\Sigma} ^{(23)}$ to denote sums over $j
\in
W_{1}, j \in\bar{W}_{1}, j \in W_{2}\cup W_{3}$ and $j
\in\bar{W}_{2}\cup\bar{W}_{3}$, respectively. Similarly,
$\Sigma^{(4)}$ denotes sum over $j \in S_{4}$, where $S_{4}= \{1,\ldots,
n\}\setminus(S_{1}\cup S_{2}\cup S_{3})$. Let $i =\sqrt{ - 1}$. Then the
following lemma is not hard to obtain using elementary facts such as
$\sin y\cos y = \frac{1}{2}\sin2y, \cos y = \frac{1}{2}(e^{iy} + e^{
- iy}),\sin y = \frac{1}{2i}(e^{iy} - e^{ - iy})$.
\begin{lemma}\label{lemaA.1}
\begin{eqnarray*}
\psi(a_{1},\ldots,a_{n})&=& \frac{1}{2^{m}i^{n_{1} +
n_{3}}}\Sigma_{1}\Sigma_{2}\Sigma_{3}M(\bar{w}_{1},\bar{w}_{2},\bar{w}_{3})
\\
&&{}\times\operatorname{exp}\biggl\{ \frac{i\pi} {2}\bigl(2\Sigma^{(1)}a_{j} - 2\bar
{\Sigma}
^{(1)}a_{j} + \Sigma^{(23)}a_{j} - \bar{\Sigma} ^{(23)}a_{j}\bigr)\biggr\},
\end{eqnarray*}
where $M(\bar{w}_{1},\bar{w}_{2},\bar{w}_{3})= ( -
1)^{\bar{w}_{1} + \bar{w}_{3}}\operatorname{exp}\{ \frac{i\pi} {4}(m -
4\bar{w}_{1}
- 2\bar{w}_{2} - 2\bar{w}_{3})\}$.
\end{lemma}

Let $g = (g_{1},\ldots,g_{n})', h = (h_{1},\ldots,h_{n})'$, where
%
%
\begin{equation}\label{equA.1}
g_{j} = u_{j} + v_{j} (\operatorname{mod} 4),\qquad h_{j} = u_{j} - v_{j}
(\operatorname{mod} 4),\qquad
1 \le j \le n.
\end{equation}

Using the same elementary facts that led to Lemma \ref{lemaA.1}, we now
note that
\begin{eqnarray}\label{equA.2}
&&\cos\biggl(\frac{\pi} {4} + \frac{\pi} {2}a'u\biggr)\cos\biggl(\frac{\pi} {4} +
\frac{\pi} {2}a'v\biggr)
\nonumber\\
&&\qquad=\frac{1}{4} \biggl[\exp\biggl\{ \frac{i\pi} {2}(1 + a'g)\biggr\}+\exp\biggl\{ \frac{i\pi}
{2}a'h\biggr\}\\
&&\qquad\quad\hspace*{11pt}{}+\exp\biggl\{ - \frac{i\pi} {2}a'h\biggr\}+\exp\biggl\{ - \frac{i\pi} {2}(1 +
a'g)\biggr\}\biggr].\nonumber
\end{eqnarray}
Throughout the rest of the \hyperref[app]{Appendix}, including the
lemmas below, we
consider $x = 0101$. Then $X = x_{1} + x_{2} + x_{3} + x_{4}= 2$, and
hence (\ref{equ2.5}), (\ref{equA.2}) and Lem\-ma~\ref{lemaA.1} yield
%
%
\begin{eqnarray}\label{equA.3}
\phi(x;a_{1},\ldots,a_{n})&=& \frac{1}{2^{({1/2})(m + 2) +
2n}i^{n_{1} +
n_{3}}}\Sigma_{1}\Sigma_{2}\Sigma_{3}M(\bar{w}_{1},\bar{w}_{2},\bar
{w}_{3})\nonumber\\[-8pt]\\[-8pt]
&&{}\times\sum_{k
= 1}^{4} \phi_{k}^{W}(x;a_{1},\ldots,a_{n}),\nonumber
\end{eqnarray}
where, with $g_{0} = 1$,
%
%
\begin{eqnarray}\label{equA.4}
&&\phi_{1}^{W}(x;a_{1},\ldots,a_{n})\nonumber\\
&&\qquad= \operatorname{exp}\biggl[\frac{i\pi}
{2}\bigl\{ g_{0} +
\Sigma^{(1)}(g_{j} + 2)a_{j} + \bar{\Sigma} ^{(1)}(g_{j} - 2)a_{j}
\\
&&\qquad\quad\hspace*{37.3pt}{}+ \Sigma^{(23)}(g_{j} + 1)a_{j} + \bar{\Sigma} ^{(23)}(g_{j} - 1)a_{j}
+ \Sigma^{(4)}g_{j}a_{j}\bigr\}\biggr],\nonumber
\end{eqnarray}
and, for $k = 2,3,4,\phi_{k}^{W}(x;a_{1},\ldots,a_{n})$ are analogous
to $\phi_{1}^{W}(x;a_{1},\ldots,a_{n})$, with $(g_{0},g_{1},\ldots
,g_{n})$ in
the latter replaced by $(0,h_{1},\ldots,h_{n}), (0, - h_{1},\ldots, - h_{n})$
and $( - g_{0}, - g_{1},\ldots, - g_{n})$, respectively. The
superscript $W$
here indicates the dependence of each
$\phi_{k}^{W}(x;a_{1},\ldots,a_{n})$ on $W_{1},W_{2},W_{3}$ via the
sums $\Sigma^{(1)},\bar{\Sigma} ^{(1)} ,\break\Sigma^{(23)} $, $\bar{\Sigma}^{(23)}$. Writing
%
%
\begin{equation}\label{equA.5}
V_{k}(x)= \frac{1}{2^{({1/2})(m + 2) + 2n}i^{n_{1} +
n_{3}}}\Sigma_{1}\Sigma_{2}\Sigma_{3}M(\bar{w}_{1},\bar{w}_{2},\bar
{w}_{3})V_{k}^{W}(x),
\end{equation}
where
%
%
\begin{equation}\label{equA.6}
V_{k}^{W}(x)=\sum_{a_{1} = 0}^{3} \cdots\sum_{a_{n} = 0}^{3}
\phi_{k}^{W}(x;a_{1},\ldots,a_{n}),\qquad 1 \le k \le4,
\end{equation}
the following lemma is immediate from (\ref{equ2.4}) and (\ref{equA.3}).
\begin{lemma}\label{lemaA.2}
$V(x)= V_{1}(x)+V_{2}(x)+V_{3}(x)+V_{4}(x)$.
\end{lemma}

Some more notation will help in presenting the subsequent lemmas. With
reference to the sets $\Delta_{ks}$ in (\ref{equ2.7}) and the $g_{j}$
and $h_{j}$
in (\ref{equA.1}), let
%
%
\begin{eqnarray}\qquad\label{equA.7}
\Delta^{(1)}&=&
\Delta_{10}\cup\Delta_{12}\cup\Delta_{30}\cup\Delta_{32},\qquad \Delta^{(2)}=
\Delta_{01}\cup\Delta_{03}\cup\Delta_{21}\cup\Delta_{23},\nonumber\\[-8pt]\\[-8pt]
\Delta&=&\Delta^{(1)}\cup\Delta^{(2)},\nonumber
\\
\label{equA.8}
\Delta_{k}^{g}&=& \{ j\dvtx1 \le j \le n, g_{j} = k\},\qquad \Delta_{k}^{h}=
\{ j\dvtx1 \le j \le n, h_{j} = k\},
\nonumber\\[-8pt]\\[-8pt]
f_{k}^{g} &=& \# \Delta_{k}^{g},\qquad f_{k}^{h} = \# \Delta_{k}^{h},\qquad 0 \le k
\le3.\nonumber
\end{eqnarray}
Then by (\ref{equ2.7}), (\ref{equ2.8}) and (\ref{equA.1}),
%
%
\begin{eqnarray}\label{equA.9}\qquad\quad
\beta_{1} &=& \# \Delta^{(1)}=f_{10} + f_{12} + f_{30} +
f_{32}=\lambda_{1} + \lambda_{3},
\nonumber\\
\beta_{2} &=& \# \Delta^{(2)}=f_{01} + f_{03} + f_{21} +
f_{23}=\lambda_{2} + \lambda_{4},
\\
\beta&=& \# \Delta=\beta_{1} + \beta_{2}=\lambda_{1} + \lambda_{2} +
\lambda_{3} + \lambda_{4},\nonumber
\\
\label{equA.10}
\Delta_{0}^{g}&=&
\Delta_{00}\cup\Delta_{13}\cup\Delta_{22}\cup\Delta_{31},\qquad \Delta_{0}^{h}=
\Delta_{00}\cup\Delta_{11}\cup\Delta_{22}\cup\Delta_{33},
\nonumber\\
\Delta_{2}^{g}&=&
\Delta_{02}\cup\Delta_{11}\cup\Delta_{20}\cup\Delta_{33},\qquad
\Delta_{2}^{h}=
\Delta_{02}\cup\Delta_{13}\cup\Delta_{20}\cup\Delta_{31},\\
\Delta_{1}^{g}&\cup&\Delta_{3}^{g} = \Delta_{1}^{h}\cup\Delta_{3}^{h}
=\Delta,\nonumber\\
\label{equA.11}
f_{2}^{g}&=&\lambda_{5} + \lambda_{7} +
\lambda_{8},\nonumber\\[-8pt]\\[-8pt]
f_{2}^{h}&=&\lambda_{6} + \lambda_{7} + \lambda_{8},\qquad
f_{1}^{g} + f_{3}^{g} = f_{1}^{h} + f_{3}^{h}=\beta.\nonumber
\end{eqnarray}
We also write
%
%
\begin{eqnarray}\label{equA.12}\hspace*{35pt}
\beta_{1}^{*}&=&f_{12} + f_{30} - f_{10} - f_{32},\qquad
\beta_{2}^{*}=f_{03} + f_{21} - f_{01} - f_{23},
\\
\label{equA.13}
\mu_{1}&=& \# \bigl(S_{3}\cap\Delta^{(1)}\bigr),\qquad \mu_{2}= \#
\bigl(S_{3}\cap\Delta^{(2)}\bigr),\\
\label{equA.14}
n(g,k)&=& \# (S_{3}\cap\Delta_{k}^{g}),\qquad n(h,k)= \#
(S_{3}\cap\Delta_{k}^{h}),\qquad 0 \le k \le3.
\end{eqnarray}
\begin{lemma}\label{lemaA.3}
\textup{(i)} $V_{1}(x)=V_{4}(x)= 0$, unless
%
%
\begin{equation}\label{equA.15}
S_{1} = \Delta_{2}^{g},\qquad S_{2}\cup S_{3}=
\Delta_{1}^{g}\cup\Delta_{3}^{g},\qquad S_{4} = \Delta_{0}^{g}.
\end{equation}

\textup{(ii)} $V_{2}(x)=V_{3}(x)= 0$, unless
%
%
\begin{equation}\label{equA.16}
S_{1} = \Delta_{2}^{h},\qquad S_{2}\cup S_{3}=
\Delta_{1}^{h}\cup\Delta_{3}^{h},\qquad S_{4} = \Delta_{0}^{h}.
\end{equation}
\end{lemma}
\begin{pf}
We prove (i). The proof of (ii) is similar. Since for
any integer~$k$,
%
%
\begin{eqnarray}\label{equA.17}
\sum_{s = 0}^{3} \exp\biggl(\frac{i\pi} {2}ks\biggr)&=& 4\qquad\mbox{if } k = 0
\operatorname{mod} 4,
\nonumber\\[-8pt]\\[-8pt]
&=& 0\qquad\mbox{otherwise},\nonumber
\end{eqnarray}
it follows from (\ref{equA.4}) and (\ref{equA.6}) that $V_{1}^{W}(x)$
vanishes, for
every $W_{1}$ $( \mbox{$\subset$} S_{1})$, $W_{2}$ $( \mbox{$\subset$} S_{2})$ and $W_{3}$ $( \mbox{$\subset$}
S_{3})$, and hence by (\ref{equA.5}) $V_{1}(x)= 0$, if either (a)
$S_{1}$ is
nonempty and $g_{j} \ne2$ for some $j \in S_{1}$, or (b) $S_{2}\cup S_{3}$
is nonempty and $g_{j} \ne1, 3$ for some $j \in S_{2}\cup S_{3}$, or (c)
$S_{4}$ is nonempty and $g_{j} \ne0$ for some $j \in S_{4}$. Thus
$V_{1}(x)= 0$, unless $S_{1} \subset\Delta_{2}^{g}, S_{2}\cup
S_{3}\subset
\Delta_{1}^{g}\cup\Delta_{3}^{g}$ and $S_{4} \subset\Delta_{0}^{g}$. These
conditions are equivalent to those in (\ref{equA.15}) because
$S_{1},S_{2},S_{3},S_{4}$ form a partition of $\{1,\ldots, n\}$ as
$\Delta_{0}^{g},\Delta_{1}^{g},\Delta_{2}^{g},\Delta_{3}^{g}$ do. The
same arguments apply to $V_{4}(x)$. Hence, (i) follows.
\end{pf}
\begin{lemma}\label{lemaA.4}
\textup{(i)} If (\ref{equA.15}) holds then
\begin{eqnarray*}
&&V_{1}(x)+V_{4}(x)\\
&&\qquad=\frac{( - 1)^{n(g,1)}}{2^{({1/2})\beta+
1}i^{n_{3}}} \biggl[\exp\biggl\{ \frac{i\pi} {4}(2 + f_{3}^{g} - f_{1}^{g})\biggr\}\\
&&\qquad\quad\hspace*{59.3pt}{}+( -
1)^{n_{3}}\exp\biggl\{ - \frac{i\pi} {4}(2 + f_{3}^{g} - f_{1}^{g})\biggr\}\biggr].
\end{eqnarray*}

\textup{(ii)} If (\ref{equA.16}) holds, then
\begin{eqnarray*}
V_{2}(x)+V_{3}(x)&=&\frac{( - 1)^{n(h,1)}}{2^{({1/2})\beta+
1}i^{n_{3}}} \biggl[\exp\biggl\{ \frac{i\pi} {4}(f_{3}^{h} - f_{1}^{h})\biggr\}\\
&&\hspace*{58.6pt}{}+( -
1)^{n_{3}}\exp\biggl\{ - \frac{i\pi} {4}(f_{3}^{h} - f_{1}^{h})\biggr\}\biggr].
\end{eqnarray*}
\end{lemma}
\begin{pf}
We prove (i), the proof of (ii) being similar. Let
(\ref{equA.15}) hold. Then by (\ref{equA.4}), noting that $g_{0} = 1$,
%
%
\begin{eqnarray}\label{equA.18}
&&\phi_{1}^{W}(x;a_{1},\ldots,a_{n})\nonumber\\
&&\qquad= \operatorname{exp}\biggl[\frac{i\pi}
{2}\bigl\{ 1 +
\Sigma^{(23)}(g_{j} + 1)a_{j} + \bar{\Sigma} ^{(23)}(g_{j} - 1)a_{j}\bigr\}\biggr]
\\
&&\qquad= \operatorname{exp}\biggl\{ \frac{i\pi} {2}(1 + 2\Sigma^{*}a_{j} +
2\Sigma^{**}a_{j})\biggr\},\nonumber
\end{eqnarray}
$\Sigma^{*}$ and $\Sigma^{**}$ being sums over $j
\in(W_{2}\cup W_{3})\cap\Delta_{1}^{g}$ and $j
\in(\bar{W}_{2}\cup\bar{W}_{3})\cap\Delta_{3}^{g}$, respectively. By
(\ref{equA.6}), (\ref{equA.17}) and (\ref{equA.18}), $V_{1}^{W}(x)= 0$
unless the ranges of
$\Sigma^{*}$ and $\Sigma^{**}$ are both empty. From this, invoking the
second equation in (\ref{equA.15}), a little reflection shows that
$V_{1}^{W}(x)=
0$ unless
%
%
\begin{eqnarray}\label{equA.19}
W_{2}&=& S_{2}\cap\Delta_{3}^{g},\qquad W_{3}=
S_{3}\cap\Delta_{3}^{g},\nonumber\\[-8pt]\\[-8pt]
\bar{W}_{2}&=& S_{2}\cap\Delta_{1}^{g},\qquad \bar{W}_{3}=
S_{3}\cap\Delta_{1}^{g}.\nonumber
\end{eqnarray}
Given $S_{2},S_{3}$, (\ref{equA.19}) determines $W_{2},W_{3} ,\bar
{W}_{2}$ and
$\bar{W}_{3}$ uniquely. Moreover, if (\ref{equA.19}) holds then $\bar
{w}_{3} =
n(g,1)$ by (\ref{equA.14}), $V_{1}^{W}(x)= 4^{n}\exp(\frac{i\pi} {2})$
by (\ref{equA.6})
and (\ref{equA.18}), and
\begin{eqnarray*}
&&
m - 4\bar{w}_{1} - 2\bar{w}_{2} - 2\bar{w}_{3}\\
&&\qquad= 2n_{1} -
4\bar{w}_{1}+w_{2} + w_{3} - \bar{w}_{2} - \bar{w}_{3}
\\
&&\qquad= 2n_{1} - 4\bar{w}_{1}+\# \{ (S_{2}\cup S_{3})\cap\Delta_{3}^{g}\}- \#
\{ (S_{2}\cup S_{3})\cap\Delta_{1}^{g}\}
\\
&&\qquad= 2n_{1} - 4\bar{w}_{1}+\# \Delta_{3}^{g}- \# \Delta_{1}^{g}= 2n_{1} -
4\bar{w}_{1}+f_{3}^{g}- f_{1}^{g}
\end{eqnarray*}
by (\ref{equA.8}), (\ref{equA.15}) and the facts that $m=2n_{1} + n_{2}
+ n_{3}, n_{k} =
w_{k} + \bar{w}_{k}$ $(k = 2, 3)$. If we summarize the above and recall
the definition of $M(\bar{w}_{1},\bar{w}_{2},\bar{w}_{3})$ from
Lem\-ma~\ref{lemaA.1}, then from (\ref{equA.5}), we get
%
%
\begin{eqnarray}\label{equA.20}
V_{1}(x)&=& \frac{4^{n}\exp({i\pi}/ {2})}{2^{({1/2})(m
+ 2) + 2n}i^{n_{1} + n_{3}}}\Sigma_{1}( - 1)^{\bar{w}_{1} +
n(g,1)}\nonumber\\[-8pt]\\[-8pt]
&&{}\times\operatorname{exp}\biggl\{ \frac{i\pi} {4}(2n_{1} - 4\bar{w}_{1} +
f_{3}^{g} -
f_{1}^{g})\biggr\}.\nonumber
\end{eqnarray}
Since $m + 2= 2n_{1}+f_{1}^{g} + f_{3}^{g}+2 = 2n_{1}+\beta+2$ by
(\ref{equA.11}) and (\ref{equA.15}), and
\[
\Sigma_{1}( - 1)^{\bar{w}_{1}}\operatorname{exp}\biggl\{ \frac{i\pi}
{4}(2n_{1} -
4\bar{w}_{1})\biggr\}=(2i)^{n_{1}}
\]
as one can verify after a little algebra, (\ref{equA.20}) yields
%
%
\begin{equation}\label{equA.21}
V_{1}(x)= \frac{( - 1)^{n(g,1)}}{2^{({1/2})\beta+
1}i^{n_{3}}} \operatorname{exp}\biggl\{ \frac{i\pi} {4}(2 + f_{3}^{g} -
f_{1}^{g})\biggr\}.
\end{equation}
Similarly, it can be shown that under (\ref{equA.15}),
%
%
\begin{equation}\label{equA.22}
V_{4}(x)= \frac{( - 1)^{n(g,3)}}{2^{({1/2})\beta+
1}i^{n_{3}}} \operatorname{exp}\biggl\{ - \frac{i\pi} {4}(2 + f_{3}^{g} -
f_{1}^{g})\biggr\}.
\end{equation}
By (\ref{equA.14}) and (\ref{equA.15}), $n_{3}=n(g,1) + n(g,3)$. Hence,
(i) follows from
(\ref{equA.21}) and (\ref{equA.22}).
\end{pf}
\begin{lemma}\label{lemaA.5}
\textup{(i)} $( - 1)^{n(h,1) - n(g,1)}=( - 1)^{\mu_{2}}$.
\textup{(ii)} If either (\ref{equA.15}) or (\ref{equA.16}) holds then
$n_{3} =
\mu_{1} + \mu_{2}$.
\end{lemma}
\begin{pf}
(i) Let $\mu_{0}= \# \{
S_{3}\cap(\Delta_{01}\cup\Delta_{23})\}$. From (\ref{equ2.7}), (\ref
{equA.1}) and (\ref{equA.8}),
\[
\Delta_{1}^{g}=
\Delta_{01}\cup\Delta_{10}\cup\Delta_{23}\cup\Delta_{32},\qquad
\Delta_{1}^{h}= \Delta_{03}\cup\Delta_{10}\cup\Delta_{21}\cup\Delta_{32}.
\]
Hence by (\ref{equA.7}), (\ref{equA.13}) and (\ref{equA.14}),
$n(h,1)-n(g,1) =\mu_{2} - 2\mu_{0}$,
and (i) follows.

(ii) If either\vspace*{1pt} (\ref{equA.15}) or (\ref{equA.16}) holds, then by (\ref
{equA.7}) and (\ref{equA.10}), $S_{3}
\subset\Delta$ $[\mbox{$=$}\Delta^{(1)}\cup\Delta^{(2)}]$. Now (ii) is immediate
from (\ref{equA.13}).
\end{pf}
\begin{lemma}\label{lemaA.6}
\textup{(i)} If $\lambda_{5} + \lambda_{6}= 0$, then
(\ref{equA.15}) and (\ref{equA.16}) become identical and both reduce to
%
%
\begin{equation}\label{equA.23}
S_{1}=\Delta_{02}\cup\Delta_{20},\qquad S_{2}\cup S_{3}=\Delta,\qquad S_{4}=
\Delta_{00}\cup\Delta_{22}.
\end{equation}

\textup{(ii)} If $\lambda_{5} + \lambda_{6}> 0$, then (\ref{equA.15})
and (\ref{equA.16}) cannot hold simultaneously.
\end{lemma}
\begin{pf}
If $\lambda_{5} + \lambda_{6}= 0$, then by (\ref{equ2.7}) and
(\ref{equ2.8}), the sets $\Delta_{11},\Delta_{13} ,\Delta_{31}$ and
$\Delta_{33}$ are empty, so that by (\ref{equA.10}) ,
$\Delta_{2}^{g}=\Delta_{2}^{h}= \Delta_{02}\cup\Delta_{20}$ and
$\Delta_{0}^{g}=\Delta_{0}^{h}=\Delta_{00}\cup\Delta_{22}$. Hence, (i)
follows recalling the last identity in (\ref{equA.10}). On the other
hand, if
$\lambda_{5} + \lambda_{6}> 0$, then at least one
of $\Delta_{11},\Delta_{13} ,\Delta_{31}$ and $\Delta_{33}$ is nonempty.
Therefore, $\Delta_{2}^{g}\ne\Delta_{2}^{h}$ and
$\Delta_{0}^{g}\ne\Delta_{0}^{h}$ and (ii) follows.
\end{pf}
\begin{lemma}\label{lemaA.7}
Let $\lambda_{5} + \lambda_{6}= 0$. Then:
\begin{longlist}
\item
$V(x)= 0$, if (\ref{equA.23}) does not hold.

\item$V(x)=( - 1)^{n(g,1)}Q_{1}(x)Q_{2}(x)$, if (\ref{equA.23})
holds, where for $k = 1, 2$,
%
%
\begin{eqnarray}\label{equA.24}
Q_{k}(x) &=& \frac{1}{2^{({1}/{2})(\beta_{k} + 1)}i^{\mu
_{k}}}\biggl[\exp\biggl\{ \frac{i\pi} {4}(1 + \beta_{k}^{*})\biggr\}\nonumber\\[-8pt]\\[-8pt]
&&\hspace*{70.2pt}{} +( - 1)^{\mu
_{k}}\exp\biggl\{ - \frac{i\pi} {4}(1 +
\beta_{k}^{*})\biggr\}\biggr].\nonumber
\end{eqnarray}
\end{longlist}
\end{lemma}
\begin{pf}
Part (i) is evident from Lemmas \ref{lemaA.2}, \ref{lemaA.3} and \ref
{lemaA.6}(i). To
prove (ii), let (\ref{equA.23}) hold and write $\alpha_{1}=\frac{\pi}
{4}(2 +
f_{3}^{g} - f_{1}^{g}), \alpha_{2}=\frac{\pi} {4}(f_{3}^{h} -
f_{1}^{h})$. Now, by Lemmas \ref{lemaA.2}, \ref{lemaA.4}, \ref{lemaA.5}
and \ref{lemaA.6}(i),
%
%
\begin{eqnarray}\label{equA.25}\hspace*{24pt}
V(x)&=&\frac{( - 1)^{n(g,1)}}{2^{({1/2})\beta+
1}i^{n_{3}}}\{e^{i\alpha_{1}}+( - 1)^{\mu_{1} + \mu_{2}}e^{ - i\alpha
_{1}}+( - 1)^{\mu_{2}}e^{i\alpha_{2}}+( - 1)^{\mu_{1}}e^{ - i\alpha
_{2}}\}\nonumber\\
&=&\frac{( - 1)^{n(g,1)}}{2^{({1/2})\beta+
1}i^{n_{3}}}\bigl\{e^{i(\alpha_{1} + \alpha_{2})/2}+( - 1)^{\mu_{1}}e^{ -
i(\alpha_{1} + \alpha_{2})/2}\bigr\}\\
&&\hspace*{0pt}{}\times\bigl\{e^{i(\alpha_{1} - \alpha_{2})/2}+(
- 1)^{\mu_{2}}e^{ - i(\alpha_{1} - \alpha_{2})/2}\bigr\}.\nonumber
\end{eqnarray}
But $\frac{1}{2}(\alpha_{1} + \alpha_{2})= \frac{\pi} {8}(2 +
f_{3}^{g} - f_{1}^{g} + f_{3}^{h} - f_{1}^{h})=\frac{\pi} {4}(1 +
\beta_{1}^{*})$, as one can verify from (\ref{equ2.7}), (\ref{equA.1}),
(\ref{equA.8}) and (\ref{equA.12}),
on simplification. Similarly, $\frac{1}{2}(\alpha_{1} -
\alpha_{2})=\frac{\pi} {4}(1 + \beta_{2}^{*})$. Hence, (ii) follows
from (\ref{equA.25}) using Lemma \ref{lemaA.5}(ii) and the fact that
$\beta= \beta_{1} +
\beta_{2}$; cf. (\ref{equA.9}).
\end{pf}
\begin{lemma}\label{lemaA.8}
\textup{(i)} If both $\beta_{1}$ and $\beta_{2}$
are odd, then the number of choices, say~$\sigma$, of
$S_{1},S_{2},S_{3}$ which meet (\ref{equA.23}) and keep both
$\mu_{1}$ and $\mu_{2}$ odd equals $2^{\beta- 2}$.

\mbox{}\phantom{i}\textup{(ii)} The number of choices of $S_{1},S_{2},S_{3}$
meeting (\ref{equA.16}) is $2^{\beta}$.

\textup{(iii)} If $\beta$ $(\mbox{$>$} 0)$ is even, then $\sigma_{0}=\sigma
_{1}=2^{\beta- 1}$,
where $\sigma_{0}$ and $\sigma_{1}$ denote the numbers of choices of
$S_{1},S_{2},S_{3}$ which meet (\ref{equA.16}) and keep $n_{3}$
even and odd, respectively.
\end{lemma}
\begin{pf}
We prove only (i). Proofs of (ii), (iii) are similar.
For a given QC design, (\ref{equA.23}) determines $S_{1}$ uniquely and fixes
$S_{2}\cup S_{3}$ at $\Delta$ [\mbox{$=$}$\Delta^{(1)}\cup\Delta^{(2)}$, by (\ref{equA.7})].
Thus, by (\ref{equA.13}), $\sigma$ equals the number of ways in which
one can
choose an odd number of elements from each of $\Delta^{(1)}$
and $\Delta^{(2)}$. Hence, $\sigma= (2^{\beta_{1} - 1})(2^{\beta_{2} -
1})=2^{\beta- 2}$, as $\beta_{1}(\mbox{$=$}\# \Delta^{(1)})$ and $\beta_{2} (\mbox{$=$}\#
\Delta^{(2)})$ are both odd and hence positive.
\end{pf}
\begin{pf*}{Proof of Theorem \ref{theo1} for $x = 0101$}
The proof is given separately
for three cases corresponding to (i1), (i2) and (i3) of Theorem \ref{theo1}.
Recall that any choice of $S_{1},S_{2},S_{3}$ which makes $V(x)$ nonzero
entails a word of type $x$, length $m + X$ and aliasing index $|V(x)|$.

\textit{Case} 1. Let $\lambda_{5} + \lambda_{6}= 0$. Then Lemma
\ref{lemaA.7} is applicable. Thus if $V(x)\ne0$ then (\ref{equA.23})
holds, so that $m +
X=2(f_{02} + f_{20})+\# \Delta+2 = l_{8}+2$, by (\ref{equ2.7})--(\ref
{equ2.9}) and (\ref{equA.9}).
Hence, each word of type $x$ has length $l_{8}+2$. It remains to show
that there are $1/(\xi_{1}^{2}\xi_{2}^{2})$ words of type $x$, each with
aliasing index $\xi_{1}\xi_{2}$. To that effect, from (\ref{equA.24})
note that
for any $S_{1},S_{2},S_{3}$ meeting (\ref{equA.23}),
%
%
\begin{equation}\label{equA.26}
|V(x)|=|Q_{1}(x)||Q_{2}(x)|,
\end{equation}
where for $k = 1, 2$,
\begin{eqnarray*}
|Q_{k}(x)|&=& 2^{ - ({1/2})(\beta_{k} - 1)}\biggl|\cos\biggl\{ \frac{\pi}
{4}(1 + \beta_{k}^{*})\biggr\} \biggr| \qquad\mbox{if } \mu_{k} \mbox{ is even}
\\
&=& 2^{ - ({1/2})(\beta_{k} - 1)}\biggl|\sin\biggl\{ \frac{\pi} {4}(1 +
\beta_{k}^{*})\biggr\} \biggr|\qquad \mbox{if } \mu_{k} \mbox{ is odd}.
\end{eqnarray*}
First, suppose both $\beta_{1}$ and $\beta_{2}$ are odd. Then by (\ref{equA.9})
and (\ref{equA.12}), $\beta_{1}^{*}$ and $\beta_{2}^{*}$ are also both
odd. For $k
= 1, 2$, if $\beta_{k}^{*}= 1$ or 5 $(\operatorname{mod} 8)$, then $|Q_{k}(x)|$ equals 0
or $2^{ - ({1/2})(\beta_{k} - 1)}$ depending on whether $\mu_{k}$ is
even or odd, while if $\beta_{k}^{*}= 3$ or 7 $(\operatorname{mod} 8)$, then the roles of
even and odd $\mu_{k}$ are switched. Hence if both $\beta_{1}^{*}$ and
$\beta_{2}^{*}$ equal 1 or 5 $(\operatorname{mod} 8)$, then by (\ref{equA.26}), a choice of
$S_{1},S_{2},S_{3}$ meeting (\ref{equA.23}) yields a nonzero $V(x)$ and hence
leads to a word if and only if $\mu_{1}$ and $\mu_{2}$ are both odd. By
(\ref{equ2.10}), (\ref{equA.9}) and Lemma \ref{lemaA.8}(i), the number
of such choices is $2^{\beta
- 2}$ $[\mbox{$=$}1/(\xi_{1}^{2}\xi_{2}^{2})]$ and, for any such choice, $|V(x)|=
2^{ - ({1/2})(\beta_{1} + \beta_{2} - 2)}$ $[\mbox{$=$}\xi_{1}\xi_{2}]$. It
is easily seen that the same holds if both $\beta_{1}^{*}$ and
$\beta_{2}^{*}$ equal 3 or 7 $(\operatorname{mod} 8)$ or if one of them equals 1 or 5
$(\operatorname{mod} 8)$ and the other 3 or 7 $(\operatorname{mod} 8)$. This settles Case 1 when
$\beta_{1}$ and $\beta_{2}$ are both odd. Similar arguments work when they
are both even or one of them is odd and the other even.

\textit{Case} 2. Let $\lambda_{5} + \lambda_{6}> 0$ and $\lambda_{1}
+ \lambda_{2} + \lambda_{3} + \lambda_{4}= 0$. By (\ref{equA.9})--(\ref
{equA.11}),
here $\beta= 0, f_{1}^{g}=f_{3}^{g} =f_{1}^{h}=f_{3}^{h} = 0$, and
$\Delta_{1}^{g}\cup\Delta_{3}^{g}$ and $\Delta_{1}^{h}\cup\Delta_{3}^{h}$
are empty sets. Hence, by (\ref{equA.14}), $n_{3}= n(g,1)= n(h,1)= 0$, for
any $S_{1},S_{2},S_{3}$ meeting (\ref{equA.15}) or (\ref{equA.16}).
Therefore, Lemma
\ref{lemaA.4}
yields $V_{1}(x)+ V_{4}(x)= 0$ under (\ref{equA.15}), and
$V_{2}(x)+V_{3}(x)= 1$
under (\ref{equA.16}). Since Lemma \ref{lemaA.6}(ii) rules out
simultaneous occurrence of
(\ref{equA.15}) and (\ref{equA.16}), Lemmas \ref{lemaA.2} and \ref
{lemaA.3} show that $V(x)$ equals 1 if
(\ref{equA.16}) holds, and 0 otherwise. Moreover, as
$\Delta_{1}^{h}\cup\Delta_{3}^{h}$ is empty, (\ref{equA.16}) determines
$S_{1},S_{2},S_{3}$ uniquely and, under (\ref{equA.16}), $m + X=
2f_{2}^{h} +
2=l_{10} + 2$, using (\ref{equ2.9}), (\ref{equA.11}). Thus, in this
case, there is one
word of type $x$, with aliasing index 1 and length $l_{10} + 2$.

\textit{Case} 3. Let $\lambda_{5} + \lambda_{6}> 0$ and $\lambda_{1}
+ \lambda_{2} + \lambda_{3} + \lambda_{4}> 0$. Again, Lemma \ref{lemaA.6}(ii)
precludes coincidence of (\ref{equA.15}) and (\ref{equA.16}). Thus, by
Lemmas \ref{lemaA.2} and \ref{lemaA.3},
words of type $x$ can arise in two mutually exclusive and exhaustive
ways: (a) from $S_{1},S_{2},S_{3}$ meeting (\ref{equA.15}), (b) from
$S_{1},S_{2},S_{3}$ meeting (\ref{equA.16}). First, consider (b). For any
$S_{1},S_{2},S_{3}$ meeting (\ref{equA.16}), by Lemmas \ref{lemaA.2},
\ref{lemaA.3} and \ref{lemaA.4}(ii),
\begin{eqnarray*}
|V(x)| &=& 2^{ - ({1/2})\beta} \biggl|\cos\biggl\{ \frac{\pi} {4}(f_{3}^{h} -
f_{1}^{h})\biggr\} \biggr|\qquad \mbox{if } n_{3} \mbox{ is even}
\\
&=& 2^{ - ({1/2})\beta} \biggl|\sin\biggl\{ \frac{\pi} {4}(f_{3}^{h} -
f_{1}^{h})\biggr\} \biggr|\qquad\mbox{if } n_{3}\mbox{ is odd.}
\end{eqnarray*}
Thus for odd $f_{3}^{h} - f_{1}^{h}$, irrespective of whether $n_{3}$ is
even or odd, $|V(x)|=2^{ - ({1/2})(\beta+ 1)}$. On the other hand,
if $f_{3}^{h} - f_{1}^{h}= 0$ or 4 $(\operatorname{mod} 8)$, then $|V(x)|$ equals $2^{ -
({1/2})\beta}$ or 0 according as whether $n_{3}$ is even or odd,
while if $f_{3}^{h} - f_{1}^{h}= 2$ or 6 $(\operatorname{mod} 8)$, then the roles of even
and odd $n_{3}$ are switched. Since $\beta$ $[\mbox{$=$}f_{1}^{h} + f_{3}^{h}$, by
(\ref{equA.11})] and $f_{3}^{h} - f_{1}^{h}$ are either both even or
both odd,
from (\ref{equ2.10}), (\ref{equA.9}) and Lemma \ref{lemaA.8}(ii),
(iii), it follows that, irrespective
of whether $f_{3}^{h} - f_{1}^{h}$ is odd or even, (b) yields
$1/(2\xi^{2})$ words of type $x$, each with aliasing index $\xi$. By
(\ref{equ2.9}), (\ref{equA.9}), (\ref{equA.11}) and (\ref{equA.16}),
each such word has length $m + X=
2f_{2}^{h}+f_{1}^{h} + f_{3}^{h}+2 =l_{10} + 2$. Similarly, (a) yields
another $1/(2\xi^{2})$ words of type $x$, each with aliasing index $\xi$
and length $l_{9} + 2$. Hence, the conclusion of Theorem \ref{theo1}
for $x = 0101$
follows in this case.
\end{pf*}
\end{appendix}

\section*{Acknowledgments}
$\!\!\!$We thank the referees for their very helpful
\mbox{suggestions}.


%

%
\printaddresses

\end{document}